\documentclass{article}

\usepackage{tikz}
	\usetikzlibrary{math}
\usepackage{authblk}
\usepackage{blindtext}
\usepackage{graphicx}

\usepackage{latexsym}
\usepackage[hidelinks]{hyperref}
\usepackage[margin=2cm]{geometry}
\usepackage{mdwlist}
\usepackage{calrsfs}

\usepackage{bm}
\usepackage{longtable}
\usepackage{pdflscape}
\usepackage{rotating}
\usepackage{algorithm}
\usepackage{algorithmic}
\usepackage[all]{xy}
\usepackage{comment}
\usepackage{url}
\usepackage{CJKutf8}

\usepackage{pgfplots, tikz} %! Packages to draw and plot

\usepackage{float}

\usepackage[cmex10]{amsmath}
\usepackage{amsthm}
\usepackage{xcolor}
\usepackage{amssymb} %For "not divide \nmid
\usepackage{amsfonts}

\usepackage[backgroundcolor=gray]{todonotes}
\usepackage{tikz}
\usetikzlibrary{chains, positioning, shapes, arrows, shapes.gates.logic.US}
\usepackage[framemethod=TikZ]{mdframed}
\newcommand{\brett}[1]{\todo[inline, backgroundcolor=green]{Brett: #1} }

\newcommand{\lucas}[1]{\todo[inline, backgroundcolor=yellow]{Lucas: #1} }

\newcommand{\old}[1]{}
\newcommand{\ignore}[1]{}
\newcommand{\subseq}[3]{#1[_{#2}^{#3}}
\newcommand{\cardinality}[1]{\#{#1}}

\begin{document}

\title{Comparing balanced % $\mathbb{Z}_v$-
sequences obtained from \\
ElGamal function to random balanced sequences}

\author[1]{Daniel Panario} 
\author[2]{Lucas Pandolfo Perin} 
\author[1]{Brett Stevens}
\affil[1]{
    School of Mathematics and Statistics,
    Carleton University, Ottawa\\
    
    ON, K1S 5B6, Canada, 
    \texttt{\{daniel,brett\}@math.carleton.ca} }
\affil[2]{
    Technology Innovation Institute, Masdar City, Abu Dhabi\\
    
    P.O. Box 9639, United Arab Emirates
    \texttt{lucas.perin@tii.ae} }
%\affil[2]{
%    Departamento de Inform\'atica e Estat\'\i stica, 
%    Universidade Federal de Santa Catarina,
%    Florian\'opolis, SC, 88040-900, Brasil
%    }

% MACROS
% \theoremstyle{definition}
\newtheorem{theorem}{Theorem}
\newtheorem{proposition}{Proposition}
\newtheorem{definition}{Definition}
\newtheorem{lemma}{Lemma}
\newtheorem{corollary}{Corollary}
\newtheorem{example}{Example}
\newtheorem*{remark}{Remark}
\newtheorem*{assumption}{Assumption}

\def\f{\mathbb{F}}
\def\fp{\mathbb{F}_{p}}
\def\fq{\mathbb{F}_{q}}
\def\R{\mathbb{R}}
\def\N{\mathbb{N}}
\def\z{\mathbb{Z}}
\def\zp{\mathbb{Z}_p}
\def\zpx{\mathbb{Z}_p^*}
\def\fpx{\mathbb{F}_p^*}
\def\rem{ \% }
\def\a{\underline{a}}
\def\b{\underline{b}}

\def\VAR{ \mbox {VAR} }
\newcommand{\Mod}[1]{\ (\mathrm{mod}\ #1)}
% make the title area
\maketitle

%%%%% TO DO %%%%%%%

\begin{abstract}
In this paper, we investigate the randomness properties of sequences in $\z_v$ derived from permutations in $\zp^*$ using the remainder function modulo $v$, where $p$ is a prime integer.  Motivated by earlier studies with a cryptographic focus we compare sequences constructed from the ElGamal function 
$x \to g^x$ for $x\in\mathbb{Z}_{>0}$ and $g$ a primitive element of 
$\zpx$, to sequences constructed from random permutations of $\zpx$.  We prove that sequences obtained from ElGamal have maximal period and 
behave similarly to random permutations with respect to the balance and run properties of Golomb's postulates 
for pseudo-random sequences. Additionally we show that they behave similarly to random permutations for the tuple balance property.  This requires some significant work determining properties of random balanced periodic sequences.  In general, for these properties and excepting 
for very unlikely events, the ElGamal sequences behave the same as random
balanced sequences.
\end{abstract}

\section{Introduction}

% Let $p$ be an odd prime number, $v \mid p-1$ and $g$ a generator of the 
% multiplicative group $\z_{p}^*$.
% In this paper we study the integer sequence defined for $1 \le i \le p-1$, by
% $$ 
% \alpha_i = (g^i \rem p) \rem v 
% $$
% where $\rem$ is the remainder function. We proceed as follows: For any permutation
% \[
%   \pi : \z_p^* \rightarrow \z_p^*,
% \]
% we are interested in properties of the sequence of the elements in $\pi_v$ where
% \[
%   \pi_v = (\pi(i) \rem v: \mbox{for } 1 \le i \le p-1 \mbox{ and } \pi(i)\in\zpx).
% \]

Applying the remainder function modulo $v$ to sequences over finite fields often yields new sequences with nice properties.  The Legendre sequence results from applying the remainder modulo $v=2$ to the sequence $\log_g(i)$ for $i \in \mathbb{F}_p = \mathbb{Z}_p$ ordered additively, where $g$ is a primitive element in $\fp$. Legendre sequences have many nice properties including low auto-correlation, high linear complexity, and favourable behaviour in the well-distribution and correlation measures \cite{GongGol2005,Mullen:2013:HFF:2555843}. Colbourn has shown that the circulant matrix formed from the same sequences modulo arbitrary $v$ is a covering array \cite{colbourn2010covering}. Tzanakis et al.~\cite{tzanakis2017} use characters over finite fields to show the existence of an infinite family of covering arrays. The elements of these covering arrays are the remainder modulo $v$ of the discrete logarithm applied to selected m-sequences, which can be viewed as the sequence of elements in $\mathbb{F}_p^{*}$ ordered multiplicatively with the trace, logarithm and remainder operators successively applied in that order.  The property required of a sequence for its circulant to yield a covering array is that every $v$-ary tuple up to some specified strength $t$ should appear at least once in some shift of every possible mask of positions with sizes up to $t$. A mask is a relative placement of positions in a circular sequence. In this sense the sequences behave as random sequences which are highly likely to have such a property when long enough.  The $t$-tuple balance property of a sequence requires that the numbers of occurrences of each tuple in the mask of $t$ consecutive positions do not differ by more than one.

The initial sequence $g^{x}$, for $1 \leq x \leq p-1$ and $g \in \zpx$ primitive, can be viewed as a permutation $g^{x}: \zpx \rightarrow \zpx$. This permutation has important applications and has been studied from multiple points of view. It is the mathematical basis for the ElGamal signature scheme \cite{elgamal}; von zur Gathen et al. have called it the {\em ElGamal function} \cite{elgamalStats}.  They point out that if this map behaves like a random permutation, then the cryptographic systems based on it are resistant to certain attacks. Motivated by this, the authors study the map $x \to g^x$ showing experiments indicating that this map behaves like a uniformly random permutation with respect to the number of cycles and  the number of $k$-cycles for a given $k$. The authors also prove, using Sidon sets, that the graph of this map is equidistributed. In their paper they explicitly call for more research of other properties for which these permutations behave like random permutations.

Drakakis et al. note that this function is also the basis of the Welch construction of Costas arrays which have very low {\em ambiguity} and are used in SONAR and RADAR frequency hopping \cite{drakakis_apn_2009,MR2723672,MR2774833}.  They refer to this function as the {\em Exponential Welch–Costas construction}. They prove \cite{drakakis_apn_2009} that this function is {\em Almost Perfect Nonlinear} which is an important property for some functions used in cryptography. They further prove that these functions are much closer to being {\em Perfect Nonlinear} than many APN functions. They point out that the SAFER cryptosystem's S-box uses an instance of the Exponential Welch–Costas construction for $p=257$. Random functions are generally nonlinear and this is another property which the map  $x \to g^x$ shares with random permutations although Drakakis et al.~\cite{MR2723672} show that the Exponential Welch–Costas construction is more nonlinear than random functions. Drakakis \cite{MR2660616} shows that the nonlinearity of random functions closely matches the nonlinearity of Costas functions for some measures, but does find two nonlinearity measures in which Costas arrays are more nonlinear than random permutations. While these results indicate that aspects of the map $x \to g^x$ can be distinguished from random, the discrepancies point to the cryptographic strength of $x \to g^x$ not its weakness.

Motivated by these investigations into the properties shared between the map $x \to g^x$ and random permutations combined with the potential power of the remainder function, we were curious about what properties would be shared between sequences derived from random permutations and the map $x \to g^x$ by using the remainder function.  We are interested in various randomness properties including balance, runs and tuple balance from Golomb's postulates, which are usually applied to sequences whose length exceeds the number of symbols \cite{GongGol2005}.

We think of the elements of $\z_p$ as equivalent classes (mod $p$), where ``mod'' is being used as a relation.  The action of applying the remainder function map $x \to x \bmod v$ is often represented with ``mod'' as a binary operator. We use $\rem$ for the remainder function to avoid confusing these two instances of mod. Using the natural cyclic additive ordering, $[1,2,\ldots,p-1]$ in $\zpx$, if $\pi: \zpx \rightarrow \zpx$ is a permutation, then $[\pi(1),\pi(2),\ldots,\pi(p-1)]$ is a periodic sequence.  Furthermore $\pi_v = [\pi(1)\rem v,\pi(2)\rem v,\ldots,\pi(p-1)\rem v]$ is a periodic sequence of length $p-1$ over $\z_v$, where $v$ can be much smaller than the length of the sequence. When $g \in \zpx$ is primitive, then $\gamma = (g^{i-1}\rem p)_{i=1}^{p-1}$ is such a permutation which can be reduced to a sequence $\gamma_v$ over $\z_v$ by taking remainders modulo $v$. This can be efficiently generated \cite{elgamalStats} and if it shares properties with random sequences over $\z_v$, then it is further evidence of the similarity of the map $x \to g^x$ to random permutations that might point towards new applications.

Before we proceed, we observe that two different representatives of the 
same $x \in \z_p$ can produce different results after applying $\rem v$. Thus in order for this transformation
to be well defined we take the underlying set of $\z_p$ to be 
$\{0, 1, \ldots , p-1\}$, that of $\z_{p}^*$ to be $\{1, 2,\ldots, p-1\}$ and 
determine $x \rem v \in \z_v$ by considering $x \in \z_p^*$ to be the corresponding 
element in $\{0, 1, \ldots, p-1\} \subset \z$.

In this paper, we are interested in randomness properties of the sequence $\gamma_v$. 
The properties we focus on are: balance, number of permutations of certain period,
$t$-tuple balance and runs. Since $\gamma$ is periodic, $\gamma_v$ is periodic as well. 
We show that if $p \equiv 1 \Mod v$, the sequences $\gamma_v$ are balanced so we compare their properties 
against random balanced periodic sequences. 

Background and 
notation are given in Section \ref{background}. In Section \ref{sec_random} we 
determine the balance, tuple and run properties for random {\em balanced} periodic 
sequences. In particular we first study the number of permutations $\pi_v$ with 
period $p-1$, showing that this number is very close to $(p-1)!$. Our main results 
of this section are proofs of the asymptotic normality for the number of times a 
$t$-tuple appears in a sequence and the number of times a symbol appears in a 
run of given length when we consider the space of random balanced sequences over 
$\z_v$ and (not necessarily minimal) period $p-1$. We also provide formulas and 
bounds for the mean and variance of these numbers. Perhaps not totally surprisingly, 
our mean and variance computations for runs are exact, while the formulas for the 
number of occurrences of tuples depend on the structure of the tuple and so we can 
only provide bounds in general.

In Section~\ref{sec_elgamal} we consider sequences which are derived from 
ElGamal permutations. We provide results for the same properties as for random 
balanced periodic sequences. We show that the balance, periodicity, tuple and 
runs of these ElGamal sequences match up nicely with the expected behaviour of 
random balanced  sequences. In Section \ref{experiments} we provide diverse 
experimental data for the number of tuples and number of runs of a symbol in 
sequences derived from ElGamal permutations. 
Finally, in Section~\ref{conclusions}, we briefly compare the studied properties
for the random permutations and the ElGamal function cases and discuss our results. 
We also comment on problems for further research including potential applications 
of this work and studies of other properties in addition to the ones presented in 
this paper.

% to add:
% 
% maguire paper and link between Sidon and Costats
% link between Sidon and APN?
% more discussion of Sidon paper?
% iterated entropy paper?
% \F_p vs \Z_p?

\section{Background and notation} \label{background}
We review some properties of sequences that are our principle interest 
in this paper. For a detailed presentation of these concepts see,
for example, \cite{GongGol2005}.

Let $\sigma$ be a permutation of $\zpx$ and $\sigma_v$ be its reduction modulo $v$.   If there exists an integer $\rho > 0$ such that 
$\sigma_v(i+\rho) = \sigma_v(i)$ for all $i$, then the sequence is \emph{periodic} 
and $\rho$ is a \emph{period of the sequence}. The smallest such $\rho$ 
is the \emph{least period of the sequence} (or simply \emph{period}).

The sequence $\sigma_v$ satisfies the \emph{balance property} if the number of times any two symbols appear in $\sigma_v$ differ by at most one. If the length of the sequence is a multiple of $v$, which is the case for our sequences,  then every
element of $\z_v$ appears the same number of times $(p-1)/v$ in a
period of the sequence.

Balance can be generalized from just a symbol to strings of length larger than one. For any sequence $s$, we denote the length $(j-i)$ subsequence $(s(i), s(i+1), \ldots, s(j-1))$ by $\subseq{s}{i}{j}$. Hence, if $s$ is length $n$ then $s = \subseq{s}{0}{n}$ and $\subseq{s}{i}{i}$ is an empty sequence.  Given a positive integer $t$, we are also interested in
counting the number of times each $t$-tuple in $\z_v^t$ appears in
the sequence.  Let $z \in \z_v^t$ and define $\lambda(z) = \cardinality{\{i:\; z=\subseq{s}{i}{i+t}\}}$. The \emph{$t$-tuple balance property} is usually expressed for m-sequences which always have length $v^n-1$ for some $n$ \cite{GongGol2005}.  Our sequences, $\sigma_v$ do not generally have these lengths and so we propose what seems to be the most natural generalization: A sequence $s$ satisfies the \emph{$t$-tuple balance property} if, for every pair $z_1,z_2 \in \z_v^t$, $|\lambda(z_1) - \lambda(z_2)| \leq 1$.

A \emph{run} of length $t$ of an element $b \in \z_v$ is a subsequence
of the form: a symbol different from $b$, followed by exactly $t$ symbols
$b$, and followed by a symbol different from $b$. The sequence satisfies 
the \emph{run property} if $(p-1)/v$ of the runs have length 1 and there are $v$ times more runs of length $i$ than runs of length $i+1$. 

We note that when $v=1$, then $\sigma_v$ is simply a sequence of 0s. When $v = p-1$, we have that $\sigma_v = \sigma$ for any permutation of $\z_p^*$. Such a sequence is balanced, with every symbol appearing once, and has period $p-1$. Only $p-1$ of the $(p-1)^t$ possible $t$-tuples appear for any $t > 1$ and every run has length one. As a consequence, for the remainder of the paper, we consider $1 < v < p-1$. For ease of reference, Table~\ref{tab:notation} gives the notation used in this paper. 
\begin{table}[ht!]
  \caption{Notation and symbols used in this paper.}
  \[\begin{array}{|rl|rl|}\hline
      p & \mbox{an odd prime}; &
      g & \mbox{a primitive element in } \z_{p}^*;\\
      v & \mbox{a positive divisor of } p-1; & 
      \rem & \mbox{the remainder function}; \\  
      a+_{n}b & (a+b) \rem n; & 
      a-_{n}b & (a-b) \rem n; \\
      a \equiv_v b &  a \equiv b \Mod{v}; &
      \pi & \mbox{a permutation in } \z_p^*; \\ 
      \pi_v & \mbox{the sequence } (\pi(1) \rem v, \ldots, \pi(p-1) \rem v); &
      |\pi_v|_a & \mbox{number of $a$'s in } \pi_v; \\
      s(i) & \mbox{the } i\mbox{-th element of the sequence } s;  & 
      \subseq{s}{i}{j} & (s(i), s(i+1), \ldots, s(j-1)); \\
      N & \mbox{maximal period } p-1; &
      \rho & \mbox{some period length at most equal to } N.
 \\ \hline
    \end{array}
  \]
\label{tab:notation}
\end{table}

The ElGamal signature scheme uses the map $x \to g^x$ to generate a signature for a message $m \in \z_{p-1}$ \cite{elgamal}. Suppose that $p$ is an $n$-bit prime number and let $g$ be a primitive element which generates $G = \z_{p-1} \cong \zpx$ with elements $[1,2,\ldots,p-1]$. The map $x \to g^x$ is a permutation so $g^x$ determines, and is determined by, $x$ uniquely.  Choose a secret global key $a \in \z_{p-1}$ uniformly at random and make $A = g^a$ public.  When signing a message $m \in \z_{p-1}$, choose a session key $k \in \z_{p-1}^{\times}$ uniformly at random.  Then, using private knowledge of $a$, compute $(K=g^k,b=k^{-1}(m-aK)) \in \z_{p-1} \times \z_{p-1}$ which is transmitted with $m$ as the signature. The signature can be verified using the public key by confirming that $g^{m} = A^KK^b$.  If for uniformly randomly chosen $(a,b) \in \z_{p-1}^2$, the triples $(A=g^a, B=g^b, C=g^{ab}) \in \z_{p-1}^3$ are indistinguishable from random triples in $\z_{p-1}$, then the ElGamal encryption scheme is resistant to public key only attacks \cite{elgamalStats}. It was this connection between random permutations and the map $x \to g^x$ which motivated some of the early work to determine properties shared between them. In acknowledgment of this history we call the map $x \to g^x$ the {\em ElGamal function}, as it is done in \cite{elgamalStats}.

% \begin{table}[ht!]
%   \caption{Notation and symbols used in this paper.}
%   \[\begin{array}{|rl|}\hline
%       p & \mbox{an odd prime}; \\
%       g & \mbox{a primitive element in } \z_{p}^*\\
%       v & \mbox{a positive divisor of } p-1; \\
%       \rem & \mbox{the remainder function}; \\  
%       a+_{n}b & (a+b) \rem n; \\
%       a-_{n}b & (a-b) \rem n; \\
%       \pi & \mbox{a permutation in } \z_p^*;  \\
%       \pi_v & \mbox{the sequence } (\pi(1) \rem v, \ldots, \pi(p-1) \rem v); \\
%       s(i) & \mbox{the } i\mbox{-th element of the sequence } s;  \\
%       N & \mbox{maximal period } p-1; \\
%       \rho & \mbox{some period length at most equal to } N. \\ \hline
%     \end{array}
%   \]
% \label{tab:notation}
% \end{table}

%%%%%%%%%%%%%%%%%%%%%%%%%%%%%%%%%%%%%%%%%%%%%%%%%%%%%%%%%%%%%%%%%%%%%%%%%%%%% 

\section{Sequences from random permutations}\label{sec_random}

In this section we take a random permutation $\pi: \zpx \rightarrow \zpx$ 
and consider the properties of the sequence $\pi_v$.
We determine the expected behaviour of interesting properties of $\pi_v$:
balance property, period length, and distribution numbers for $t$-tuples 
and runs over all random balanced permutations $\pi: \zpx \rightarrow \zpx$.

%%%%%%%%%%%%%%%%%%%%%%%%%%%%%%%%%%%%%%%%%%%%%%%%%%%%%%%%%%%%%%%%%%%%%%%%%%%%% 
\subsection{Balance property and number of permutations $\pi_v$ with period $\rho$}\label{sec_balance_random}
Because a permutation of $\zpx$ has length $p-1$ and the set of values in the 
interval $[1,p-1]$ which are congruent to $i \Mod v$ has cardinality 
$\lceil (p-1 - ((i-1) \bmod v))/v\rceil$, the balance property of $\pi_v$ is 
readily established.
\begin{proposition} \label{balancedproperty}
  Let $\pi$ be a permutation in $\zpx$, then $\pi_v$ is balanced. Moreover $|\pi_v|_a = |\pi_v|_b$ for all symbols $a, b \in \z_v$ if and only if $v \mid p-1$; we call this {\em exactly balanced}.
\end{proposition}

We are most interested in exactly balanced sequences $\pi_v$ and so we focus on the case that $p \equiv 1 \Mod v$.  However some of our results generalize nicely to any $1<v<p$. We explicitly note how these more general cases behave and some results are stated in their full generality.

% %%%%%%%%%%%%%%%%%%%%%%%%%%%%%%%%%%%%%%%%%%%%%%%%%%%%%%%%%%%%%%%%%%%%%%%%%%%%% 
% \subsection{Number of permutations $\pi_v$ with period $\rho$}
In this paper we use the term {\em period} to refer to the smallest possible 
period of the sequence. When $p \not\equiv 1 \Mod v$, then $\pi_v$ has period 
$p-1$ for any $\pi: \zpx \rightarrow \zpx$, as we show next.
\begin{lemma} \label{lemma_near_balanced_period}
If $p \equiv \alpha \Mod v$ for $\alpha \neq 1$, then $\pi_v$ has period $p-1$ for any $\pi: \zpx \rightarrow \zpx$.
\end{lemma}
\begin{proof}
    If $p \equiv \alpha \Mod v$ where $\alpha \neq 1$, then
  \[
    |\pi_v|_a = \begin{cases} \lceil (p-1)/v \rceil &  0 \leq a < \alpha -1, \\
      \lfloor (p-1)/v \rfloor & \mbox{otherwise}.
    \end{cases}
    \]
    Suppose that $\pi_v$ has period $N$. Then the number of times any symbol occurs in $\pi_v$ is equal to the number of times it appears in any period multiplied by the number of periods.  Thus since $\alpha \neq 1$,
    \[
  1 =     |\pi_v|_0 - |\pi_v|_{v-1} = \frac{p-1}{N} \left(|\subseq{\pi_v}{0}{N}|_0 - |\subseq{\pi_v}{0}{N}|_{v-1} \right) \geq \frac{p-1}{N} \geq 1,
\]
because $|\subseq{\pi_v}{0}{N}|_0 - |\subseq{\pi_v}{0}{N}|_{v-1}$ is a positive integer.
We conclude that $N = p-1$.
\end{proof}
We have shown that $\pi_v$ can only have period smaller than $p-1$ when $p \equiv 1 \Mod v$.  Next we show that even in this case smaller periods are exceptionally rare.

\begin{lemma}\label{many_to_one}
If $\sigma$ is a balanced sequence over $\z_v$ of length $p-1$, then there 
are exactly $(((p-1)/v)!)^v$ permutations $\pi$ of $\zpx$ such that 
$\pi_v = \sigma$.
\end{lemma}
\begin{proof}
  We have that $v \mid p-1$, since $\sigma$ is balanced. Let $j$ be a value 
  in $\z_v$. Then, the elements in the $(p-1)/v$ positions of $\sigma$ that are equal 
  to $j$ can be lifted to elements in $\zpx$ in $((p-1)/v)!$ ways. This happens for
  every distinct element $j$ in $\z_v$. Hence, the total number of permutations
  $\pi$ in $\zpx$ where $\pi_v = \sigma$ is 
  $$ \left(\left(\frac{p-1}{v}\right)!\right)^v.$$
\end{proof}
Let $s(v,\rho)$ be the number of balanced sequences $\pi_v$ over $\z_{v}$
satisfying $\pi_v(i+\rho) = \pi_v(i)$ for all $0 \leq i < \rho$. Let $t(v,\rho)$ be the number of balanced sequences $\pi_v$ over $\z_{v}$
with period $\rho$.

\begin{proposition}
For $s(v,\rho)$ and $t(v,\rho)$ as defined above, we have
  \begin{equation*}
    s(v,\rho) = \binom{\rho}{\rho/v, \rho/v, \ldots, \rho/v}, \quad    
    t(v,\rho) = \sum_{d | \frac{\rho}{v}} \mu(d) s(v,\rho/d) 
  = \sum_{d| \frac{\rho}{v}} \mu(d) \binom{\rho/d}{\rho/vd, \rho/vd, \ldots, \rho/vd},
  \end{equation*}
  where $\mu$ is the M\"{o}bius function.
\end{proposition}
\begin{proof}
  For $s(v,\rho)$, in each $\rho$ consecutive entries 
  of a balanced sequence with $\pi_v(i+\rho) = \pi_v(i)$ there must be 
  exactly $\rho/v$ occurrences of each symbol in $\z_v$.  Thus $\rho$ must 
  be a multiple of $v$ and these sequences are enumerated by multinomial coefficients.

  Similarly for $t(v,\rho)$, for every balanced sequence over
  $\z_{v}$ with period $\rho$ we must have that $\rho$ is a multiple 
  of $v$. We have that $s(v,\rho) = \sum_{d | \rho} t(v,d)$ and therefore the 
  M\"{o}bius Inversion Formula gives the result.
\end{proof}

\begin{corollary} \label{nummaxperiod}
Let $p \equiv 1 \Mod v$ be prime and $\pi$ be a permutation in $\zpx$. The number of permutations
  $\pi_v$ with period $p-1$ is
  \[
    T = t(v, p-1) \left( \left(\frac{p-1}{v}\right) ! \right)^v.
  \]
\end{corollary}

\begin{theorem}
  Let $Q$ be the set of all prime divisors of $k = (p-1)/v$, and $q'$ be the 
  minimal element of $Q$. Then, the number of permutations $\pi_v$ with period 
  $p-1$ satisfies 
  \begin{equation} \label{T_bound}
    (p-1)! \left( 1 - \sum_{q \in Q} \frac{(\frac{p-1}{q})! (k!)^v}
      {(p-1)! \left(\frac{k}{q}\right)!\ldots \left(\frac{k}{q}\right)!}\right)
    \le T \le
    (p-1)!\left(1 - \frac{(\frac{p-1}{q'})! (k!)^v}
      {(p-1)! \left(\frac{k}{q'}\right)!\ldots \left(\frac{k}{q'}\right)!} \right ).
  \end{equation}
  As a consequence, for every $\epsilon$ there exists an $n_{\epsilon}$ so that for 
  all $p \geq n_{\epsilon}$,
  \begin{equation} \label{T_bound_epsilon}
    (p-1)!(1-\epsilon) \leq T \leq (p-1)!.
  \end{equation}
\end{theorem}
\begin{proof}
  Since $\mu(d)$ is $-1$ for $d$ prime, we claim that considering only $d=1$ and the prime divisors of $\rho=p-1$ gives a lower bound for $T$.  To see this we first note that for $v \geq 2$ an induction proves
  \[
    \sum_{i = 0}^{k} \binom{vi}{i,i,\ldots,i} \leq \binom{v(k+1)}{k+1,k+1,\ldots,k+1}.
  \]
  A consequence of this is that for any prime $q | (\rho/v) $
  \[
    \sum_{q|d| \frac{\rho}{v} \atop q < d}  \binom{\rho/d}{\rho/vd, \rho/vd, \ldots, \rho/vd} \leq \binom{\rho/q}{\rho/vq, \rho/vq, \ldots, \rho/vq}
  \]
  which is enough to prove our claim.

  Hence, letting $k = (p-1)/v$ and 
  $q \in Q = \{ q: q \mid k \mbox{ and }q \mbox{ is prime} \}$ we have 

  \begin{align*}
    T & \ge \left( \binom{p-1}{k, \ldots, k} + \sum_{q \in Q} \mu(q)\binom{(p-1)/q}{k/q, \ldots k/q} \right) (k!)^v\\
      &= (p-1)!\left (1 - \sum_{q \in Q} \frac{(\frac{p-1}{q})! (k!)^v}{(p-1)! \left(\frac{k}{q}\right)!\ldots \left(\frac{k}{q}\right)!} \right).
  \end{align*}

  Similarly, for the upper bound, we proceed with the same rationale but consider only 
  the smallest prime divisor $q'$ of $k$ to obtain the bound in Equation~(\ref{T_bound}).

  For the asymptotic result, from \cite{MR69328} we have, for all $n \geq 1$,
  \[
    \sqrt{2\pi} n^{n+\frac{1}{2}} e^{-n} e^{\frac{1}{12n +1}} \leq n! \leq \sqrt{2\pi} n^{n+\frac{1}{2}} e^{-n} e^{\frac{1}{12n}}.
  \]
  On the left hand side we use that, for all $n \geq 1$, $e^{1/(12n+1)} \ge 1$. 
  On the right hand side we first observe that $e^{(12n-1)/12n} > \sqrt{2\pi}$ 
  for all $n \geq 2$, and then we note that the following upper bound also holds 
  for $n=1$. Hence, for all $n \geq 1$, we have the bounds 
  \[
    \sqrt{2\pi} n^{n+\frac{1}{2}} e^{-n} \leq n! \leq e n^{n+\frac{1}{2}} e^{-n}.
  \]
  From these bounds, and recalling that $v=(p-1)/k$, we obtain that
  \[
     \left(\frac{\sqrt{2\pi}}{e}\right)^{v+1} \frac{q^{\frac{v-1}{2}}}{v^{(p-1)(1-1/q)}} 
\leq \frac{(\frac{p-1}{q})!(k!)^v} {(p-1)! \left(\left(\frac{k}{q}\right)!\right)^v} 
\leq \left(\frac{e}{\sqrt{2\pi}}\right )^{v+1} \frac{q^{\frac{v-1}{2}}}{v^{(p-1)(1-1/q)}}.
\]

Recalling that $q'$ is the smallest element in $Q$, we substitute 
these bounds into Equation~(\ref{T_bound}) getting
\[
  (p-1)! \left( 1 - \sum_{q \in Q} \left(\frac{e}{\sqrt{2\pi}}\right )^{v+1} \frac{q^{\frac{v-1}{2}}}{v^{(p-1)(1-1/q)}} \right)
    \le T \le
    (p-1)!\left(1 -   \left(\frac{\sqrt{2\pi}}{e}\right)^{v+1} \frac{{q'}^{\frac{v-1}{2}}}{v^{(p-1)(1-1/q')}}  \right ).
  \]
  We note that $2 \leq q \leq k$ for all $q \in Q$  and that 
  $-c (q^{(v-1)/2})/(v^{(p-1)(1-1/q)})$ has positive derivative on the interval 
  $[2,k]$ for any positive $c$. Thus,
\[ (p-1)! \left( 1 - \cardinality{Q} \left(\frac{e}{\sqrt{2\pi}}\right )^{v+1}
\frac{{q'}^{\frac{v-1}{2}}}{v^{(p-1)(1-1/q')}} \right) \le T \le
(p-1)!\left(1 - \left(\frac{\sqrt{2\pi}}{e}\right)^{v+1}
\frac{{q'}^{\frac{v-1}{2}}}{v^{(p-1)(1-1/q')}} \right ).
  \]\
Using the fact that $\cardinality{Q} \leq (p-1)/v$, we obtain
\[ (p-1)! \left( 1 -  \left(\frac{e}{\sqrt{2\pi}}\right )^{v+1}
\frac{(p-1){q'}^{\frac{v-1}{2}}}{vv^{(p-1)(1-1/q')}} \right) \le T \le
(p-1)!\left(1 - \left(\frac{\sqrt{2\pi}}{e}\right)^{v+1}
\frac{{q'}^{\frac{v-1}{2}}}{v^{(p-1)(1-1/q')}} \right ).
  \]
Finally we use $2 \leq q' \leq (p-1)/v$ to  obtain the simplified lower and upper bounds 
  \[
    (p-1)!\left(1 - \left(\frac{e}{\sqrt{2\pi}}\right)^{v+1}
          \frac{(p-1)2^{\frac{v-1}{2}}}{v^{(p+1)/2}} \right)
    \le T \leq
    (p-1)!\left(1 - \left(\frac{\sqrt{2\pi}}{e}\right)^{v+1}
          \frac{(p-1)^{\frac{v-1}{2}}}{v^{(2p-v-3)/2}}\right).
  \]

  The intricate fraction is polynomial in $p$ in the numerator and exponential in 
  $p$ in the denominator so, for a fixed $v$ and $p$ even moderately large, $T$ 
  is very close to $(p-1)!$.
\end{proof}

When $p = vq+1$ for $p$ and $q$ prime, we can compute the exact values. We have 
that $\rho$ must be both a multiple of $v$ and a divisor of $p-1 = vq$ and so 
its only two values are $v$ and $vq$:
\begin{align*}
  s(v,v) = v!, \qquad s(v,vq) = \binom{vq}{q,\ldots,q},
\end{align*}
and
\begin{align*}
  t(v,v) = v!, \qquad t(v,vq) = \binom{vq}{q,\ldots,q} - v!. 
\end{align*}
The number of permutations of $\zpx$, out of $(p-1)!$ which give a sequence 
with period shorter than $p-1 = vq$ is
\begin{equation} \label{a048617}
  v! \left( \frac{vq}{v} !\right )^v =  v!(q!)^v.
\end{equation}
Thus, the chance that a random permutation of $\zpx$ gives a sequence with 
period shorter than $p-1 = vq$ is
\[
  \frac{v!(q!)^v}{(vq)!} = \frac{v!}{\binom{vq}{q,\ldots,q}},
\]
which is generally very small. When $v=2$ these primes $p$ are known as 
{\em safe primes} (while the prime number $q$ such that $p=2q+1$ is known as 
a Sophie Germain prime). When $v=2$ and $p$ is only as high as $83$, the 
probability that a random sequence has period less than the full period is 
already extremely small: $t(v,v)/(p-1)! = 4.708 \times 10^{-24}$. For $v=8$ 
and $p=137$, we calculate it as $t(v,v)/(p-1)! = 2.82 \times 10^{-112}$. 

% Table~\ref{tab:occurencesPeriodic} shows these probabilities for $v \leq 8$ 
% and $p < 200$. \daniel{We need to shorten this table with few interesting or
% representative values only.}
% \begin{table}[ht!]
%   \caption{Table of number of permutations in $\f_p^*$ with period $N=2$ and
%     the ratio with the total number of permutations in $\fp^*$.}
%   \[\begin{array}{rrl|rrl}\hline
%       v & p & t(v,v)/(p-1)! & v & p & t(v,v)/(p-1)! \\ \hline\hline
%       2 & 5 & 3.333 \times 10^{-1} &
%       2 & 7 & 1.00 \times 10^{-1} \\
%       2 & 11 & 7.937 \times 10^{-3} &
%       2 & 23 & 2.835 \times 10^{-6} \\
%       2 & 47 & 2.429 \times 10^{-13} &
%       2 & 59 & 6.652 \times 10^{-17} \\
%       2 & 83 & 4.708 \times 10^{-24} &
%       3 & 7 & 6.667 \times 10^{-2}\\
%       4 & 13 & 6.494 \times 10^{-5} &
%       4 & 29 & 5.079 \times 10^{-14} \\
%       4 & 53 & 4.474 \times 10^{-28} &
%       4 & 149 & 3.369 \times 10^{-85} \\
%       4 & 173 & 1.498 \times 10^{-99} &
%       5 & 11 & 1.058 \times 10^{-3}\\
%       6 & 13 & 9.620 \times 10^{-5} &
%       6 & 19 & 5.247 \times 10^{-9} \\
%       6 & 31 & 8.105 \times 10^{-18} &
%       6 & 43 & 8.399 \times 10^{-27} \\
%       6 & 67 & 5.350 \times 10^{-45} &
%       6 & 79 & 3.707 \times 10^{-54} \\
%       6 & 103 & 1.516 \times 10^{-72} &
%       6 & 139 & 3.107 \times 10^{-100} \\
%       8 & 17 & 4.933 \times 10^{-7} &
%       8 & 41& 2.125 \times 10^{-27} \\
%       8 & 89 & 1.401 \times 10^{-69} &
%       8 & 137 & 2.820 \times 10^{-112}
%     \end{array}
%   \]
% \label{tab:occurencesPeriodic}
% \end{table}

We note that the value in Equation~(\ref{a048617}) is the order of the 
automorphism group of the complete multipartite graph $K_{q,\ldots,q}$ with 
$v$ parts of size $q$, and for $v=2$ is sequence (A048617) at the 
Online Encyclopedia of Integer Sequences\footnote{http://oeis.org/A048617}.

\subsection{$t$-tuple balance and run properties of random balanced sequences} \label{random_tuples}

We now want to understand the tuple balance and run properties of $\pi_v$. Since 
Lemma~\ref{many_to_one} shows that $(((p-1)/v)!)^v$ permutations of $\zpx$ yield 
the same sequence over $\z_v$, in the remainder of this section we focus simply 
on balanced sequences $\sigma$ of length $n$ over $\z_v$.  For $z \in \z_v^t$ a 
$t$-tuple over $\z_v$ and $\sigma \in B(v,n) = \{\sigma \in \z_v^n: \sigma 
\mbox{ is balanced}\}$, let
\[
  \lambda(z) = \cardinality{\{i \in [1,p-1]: \sigma(i+_{n} \iota) = z(\iota),\; 0 \leq \iota < t\}},
\]
where we recall that $i+_{n} \iota = (i+\iota) \rem n$. 
This counts the number of times a fixed $t$-tuple appears in $\sigma$. For 
$b \in \z_v$, $t \in \mathbb{N}$ and $\sigma \in B(v,n)$, to count the runs 
of symbol $b$ and length $t$, let
\[
  \rho(b,t) = \cardinality{ \{i \in [1,p-1]: \sigma(i-_{n} 1) \neq b,\; \sigma(i+_n t) \neq b,\;  
  \sigma(i+_n\iota) = b,\; 0 \leq \iota < t\}}.
\]

\subsubsection{Asymptotic normality} \label{asymp_normal}

We are interested in the distribution of values of $\lambda(z)$ and $\rho(b,t)$ 
as $\sigma$ ranges over all balanced sequences in $B(v,n)$. We show that for 
$n=vl$, as $l \rightarrow \infty$, these values are distributed normally; 
we also compute and bound the mean and variance. We follow closely the framework 
of Bender, Richmond and Williamson \cite{MR721368}; readers interested in the 
full details may find helpful to have reference~\cite{MR721368} at hand. 
Using their Theorem~1 we establish central and local limit theorems for a higher 
dimensional distribution. Then, we restrict to the case of interest using 
conditional density to guarantee that the distribution is still normal 
\cite[Proposition~3.13]{10.2307/20461449}. 
Finally, in Section \ref{mean_variance}, we directly compute and bound the 
mean and variance.

The techniques required in our proofs of this section introduce generating 
functions and use the \emph{transfer matrix method}. For a more detailed 
introduction of the transfer matrix method and its many counting applications
see \cite[Section V.6]{FlajSed2009}. 
We use the following notation: for $x = (x_0,x_1,\ldots,x_{v-2},u)$ and 
$k \in \mathbb{N}^v$ given by $k = (k_0,k_1,\ldots,k_{v-1})$, we write 
$x^k = (\prod_{0 \leq b < v-1}x_b^{k_b} )u^{k_{v-1}}$.  

We are interested in counting circular sequences in $B(v,n)$. Hence, for 
$Z \subset \z_v^t$, we define the $v^t \times v^t$ transfer matrix $T$ with 
rows and columns indexed by $\mathbb{Z}_v^t$:
\[
  T_{z',z''} = \begin{cases}
    u x_{z'(0)}  & \mbox{if } z' \in Z,\; z^{'}(0) \neq v-1,\; z'(i+1) = z''(i)\;\; 0 \leq i < t-1; \\
    u   & \mbox{if } z' \in Z,\; z^{'}(0) = v-1,\; z'(i+1) = z''(i)\;\; 0 \leq i < t-1; \\
    x_{z'(0)}  & \mbox{if } z' \not\in Z,\; z^{'}(0) \neq v-1,\; z'(i+1) = z''(i)\;\; 0 \leq i < t-1; \\
    1  & \mbox{if } z' \not\in Z,\; z^{'}(0) = v-1,\; z'(i+1) = z''(i)\;\; 0 \leq i < t-1; \\
    0 & \mbox{otherwise}.
  \end{cases}
\]
To count circular sequences, we also define the $v^t \times v^t$ matrix $C$ 
with rows and columns indexed by $\mathbb{Z}_v^t$:
\[
  C_{z',z''} = \begin{cases}
    1  & \mbox{if } z' = z'', \\
    0 & \mbox{otherwise.}
  \end{cases}
\]
For $k = (k_0, \ldots, k_{v-1})$, the generating function counting the 
numbers $a_n(k)$ of circular sequences in $\z_v^n$ which have $k_b$ 
occurrences of symbol $b$ for $0 \leq b < v-1$, $n-\sum_{b=0}^{v-2} k_b$ 
occurrences of symbol $v-1$, and in which $k_{v-1}$ is the number of times 
strings from $Z$ appear, is given by
\[
  \sum_{k \in \mathbb{N}^t} a_n(k)x^k = \sum_{z',z'' \in \z_v^t} C_{z',z''} T^{n}_{z',z''}.
\]
We observe that the parameter $u$ in the generating function marks the latter
condition, that is, that the number of times a string from $Z$ appears is 
exactly $k_{v-1}$. We also remark that $C_{z',z''} T^{n}_{z',z''}$ is an 
entrywise multiplication, not matrix multiplication, that given the shape of 
$C_{z',z''}$ leads to the counting of circular sequences of length $n$. 

Bender et al.~\cite{MR721368} determine some sufficient conditions for this 
type of generating function to have its values distributed normally. We do 
not need the full strength of their theorem, so to simplify the statement 
we first let $D$ be the directed graph with vertex set $\z_v^t$ and containing 
an arc $(z',z'')$ if $T_{z',z''} \neq 0$. A de Bruijn cycle is a Hamilton 
cycle in $D$ and these exist for every $v$ and $n$ \cite{MR47311}. This shows 
that $D$ consists of a single strongly connected component. Under these 
conditions, Bender et al.'s definitions can be simplified as follows. 
\begin{definition}\cite{MR721368}
% Suppose that $r_l > 0$ for all $l$. 
With the above notation, the sequence of numbers $a_n(k)$ is 
{\em admissible at $r$ for $n \equiv n_0 \Mod d$}, where every entry of $r$ is positive, if
\begin{enumerate}
\item the entries in $C$ and $T$ have Laurent series expansions about 0 
having nonnegative coefficients and converging in a neighborhood of $r$; 
\item the greatest common divisor of the lengths of the directed cycles 
of $D$ is $d$.
\end{enumerate}
\end{definition}
\begin{definition}\cite{MR721368} \label{bender_def_2}
Suppose the sequence of numbers $a_n(k)$ is admissible. For $z',z'' \in \z_v^t$ 
let $\mathcal{A}_{z',z''}^{(s)}$ be the additive Abelian group generated 
by vectors of the form $k_1 - k_2$, where $k_1$ and $k_2$ are exponents
of terms in the $(z',z'')$ entry of $T^s$. If the entry is zero, set
$\mathcal{A}_{z',z''}^{(s)} = \emptyset$. Let $\Lambda$ be the lattice  
\[
 \Lambda =  \cup_{s,z',z''} \mathcal{A}_{z',z''}^{(s)}.
\]
\end{definition}
We have simplified the statements of these two definitions given 
in~\cite{MR721368} because in our 
case $D$ is strongly connected. Bender et al. then prove the following 
important result.
\begin{theorem}[Theorem 1 of \cite{MR721368}] \label{bender_thm}
Suppose $a_n(k)$ is admissible at 1 for $n \equiv n_0 \Mod d$ and that 
$\Lambda$ is d-dimensional. Then $a_n(k)$ satisfies a central limit theorem 
for $n \equiv n_0 \Mod d$ with means and covariance matrix asymptotically 
proportional to $n$.  Let $q$ be such that $qc \in \Lambda$ for all 
$c \in \z^v$.  Then $a_n(k)$ satisfies a local limit theorem modulo 
$\Lambda$ for $n \equiv n_0 \Mod {dq}$
\end{theorem}

To invoke Theorem~\ref{bender_thm}, we require that the numbers $a_n(k)$ 
be admissible. As pointed out in \cite{MR721368}, the fact that the entries 
of $T$ and $C$ are polynomials in the $x_i$ and $u$ implies that the 
first condition of admissibility is satisfied. To show that the greatest 
common divisor $d$ of cycle lengths in $D$ is one, consider the strings 
$i^tj^si^t$. These  correspond to cycles starting and finishing at vertex 
$i^t$ with length $2t+s$. By choosing two consecutive values of $s$ 
we conclude that $d=1$ and thus $a_n(k)$ is admissible with $d=1$. 

All that remains is to compute the lattice $\Lambda$ given in 
Definition~\ref{bender_def_2}. We show next that $\Lambda = \z^v$ 
for a large set of possible $Z$.
\begin{lemma}\label{lattice_lemma}
If there exists $w \in \z_v^{t-2}$, $a \neq b \in \z_v$ and $c \neq d \in \z_v$ 
such that the cardinality of
  \[
    \{awc,awd, bwc, bwd\} \cap Z
  \]
  is odd, then $\Lambda = \z^v$. 
\end{lemma}
\begin{proof}
Suppose $w, a, b, c, d$ satisfy the required properties; let $z \in \z_v^t$ 
and $x \in \z_v$ be arbitrary and fixed. Since counting the occurrences of 
tuples from $Z$ in a word uniquely determines the number of tuples from 
$Z^c$, the complement of $Z$, in the word, we assume that only one of 
$\{awc,awd, bwc, bwd\}$ is in $Z$ and without loss of generality it is $awc$. 

For any $v-1 \neq y \in \z_v$, let $e_y$ denote the elementary unit vector of 
length $v$ with a 1 in the $y$th position and let $e_u$ denote the elementary 
unit vector of length $v$ with a 1 in the last position. Since there is no 
variable for symbol $v-1$, $e_{v-1}$ shall denote the length $v$ vector of 0s.
In $T_{z,wcx}^{t+1}$ are the terms $p_zx_au^e$ and $p_zx_bu^f$ where $p_z$ 
depends only on $z$ and contains no $u$, and $e$ and $f$ depend on $z$, $w$ 
and $a$, $b$ and $c$. 
Thus $e_a-e_b+(e-f)e_u \in \mathcal{A}_{z,wcx}^{(t+1)} < \Lambda$. In 
$T_{z,wdx}^{t+1}$ are the terms $p_zx_au^{e-1}$ and $p_zx_bu^{f}$. Thus 
$e_a-e_b +(e-f-1)e_u\in \mathcal{A}_{z,wdx}^{(t+1)} < \Lambda$ and hence, 
$e_u \in \Lambda$.

For any $i,j \in \z_v$, there are terms $p_zx_iu^g$ and $p_zx_ju^h$ in 
$T_{z,wcx}^{t+1}$. Thus, $e_i-e_j+(g-h)e_u \in \Lambda$ and since 
$e_u \in \Lambda$ we have $e_i - e_j \in \Lambda$. By setting $j=v-1$, we have 
$e_i \in \Lambda$ for all $i \in \z_v, i \neq v-1$. This completes the proof.
\end{proof}

\begin{corollary}
If $\cardinality{Z} = 1$, then $\Lambda = \z^v$.
\end{corollary}

\begin{proof}
Suppose $Z = \{x_0x_1\cdots x_{t-1} \}$. Let $a = x_0$, $c = x_{t-1}$, 
$w = x_1x_2 \cdots x_{t-2}$, choose $b \neq a$ and $d \neq c$, and apply Lemma~\ref{lattice_lemma}.
\end{proof}

\begin{corollary}
Let $x \in \z_v$ and suppose $Z = \{z \in \z_v^{t+2}:\; z(0) \neq x, 
z(t+1) \neq x,\; z(\iota) = x, 1 \leq \iota < t+1\}$. Then $\Lambda = \z^v$.
\end{corollary}

\begin{proof}
Let $a\neq x,c \neq x$, $b=d=x$ and $w = xx\cdots x$ and apply Lemma~\ref{lattice_lemma}.
\end{proof}

Although not every possible $Z$ is dealt with by Lemma~\ref{lattice_lemma}, the overwhelming majority of possible sets $Z$ is covered by the lemma, as the next theorem shows.
\begin{theorem}
	The number of $Z \subset \z_v^t$ which satisfy Lemma~\ref{lattice_lemma} is $2^{(v-1)^{2}v^{t-2}}$.
\end{theorem}
\begin{proof}
We count the number of sets $Z$ which fail to satisfy Lemma~\ref{lattice_lemma}. Such 
a $Z$ has the property that for every $w \in \z_v^{t-2}$, every $a\neq b \in \z_v$ 
and every $c \neq d \in \z_v$, an even number of $\{awc,awd,bwc,bwd\}$ are in $Z$. 
For each $w$, let $Z_w = \{xwy \in Z: x,y \in \z_v\}$. As $Z$ varies over all possible subsets of $\z_v^{t}$, we observe that the number $m$ of 
$Z_w$ which do not satisfy the lemma does not depend on $w$.  Furthermore, if both $Z_w$ and 
$Z'_w$ fail to satisfy the lemma, then $(Z\setminus Z_w)\cup Z'_w$ satisfies 
Lemma~\ref{lattice_lemma} if and only if $Z$ does. We conclude that the number of 
$Z$ which do not satisfy Lemma~\ref{lattice_lemma} is $m^{v^{t-2}}$.

For a particular $w \in \z_v^{t-2}$, if $Z_w$ and ${Z'}_w$ both fail to satisfy the 
lemma, then so does their symmetric difference $Z_w \triangle {Z'}_w$. This means 
that the incidence vectors of sets $Z_w$ which fail to satisfy the lemma form a 
code of length $v^2$. It is easy to check that for any fixed $a$ the incidence 
vectors of both
  \[
    P_{aw*} = \{awx: x \in \z_v\}, 
  \]
  and
  \[
    S_{*wa} = \{xwa: x \in \z_v\}, 
  \]
are in the code. For any vector $v$ in the code we add a linear combination of 
vectors from $\{S_{*wa}:a \in \z_v\}$ to obtain a vector that has a 0 in all the positions 
given by $0wx$, $x \in \z_v$. Suppose that there is a 1 in position $awx$ for 
$a \neq 0$. The vector $0wx$ is not in the set and for every $y$, neither is $0wy$.  Since the set represented by this vector fails to satisfy the lemma, $awy$ must be in this set.  Thus there must be a 1 in position $awy$ for all $y \in \z$. 
Hence, for each $a \neq 0$, the values in positions $awx$ are the same for all $x$. Therefore, 
the incidence vectors of $P_{aw*}$ and $S_{*wa}$ span the space and the code has 
dimension $2v-1$ and $m = 2^{2v-1}$. Thus, the number of $Z$ which fail to satisfy the lemma is
  \[
    2^{(2v-1)v^{t-2}}.
  \]
Correspondingly the number of $Z$ which do satisfy the lemma is
  \[
    2^{(v^2-2v+1)v^{t-2}}.
  \]
\end{proof}

In either scenario, tuples or runs, all the conditions of Bender et al.'s Theorem~1 
\cite{MR721368} are met. Moreover $d=1$ and the lattice $\Lambda = \z^v$.  Thus we can conclude that the distribution of $a_n(k)$ is asymptotically joint normal for every $n$.
\begin{theorem}\label{thm_normal_random}
Let $z \in \z_v^t$ and $t(\kappa)$ be the number of balanced circular sequences 
of length $n$ over $\z_v$ for which $\lambda(z) = \kappa$. 
There exists a $m_{\lambda},b_{\lambda},c_{\lambda} \in \R$ such that
  \[
    \sup_{\kappa} \left |\frac{\sqrt{2\pi b_{\lambda}}t(\kappa)}{\binom{vl}{l,l,\ldots,l}}-c_{\lambda} e^{(\kappa-m_{\lambda})^2/b_{\lambda}} \right| = o(1). 
  \]
Let $b \in \z_v$, $t \in \N$ and $r(\kappa)$ be the number of balanced circular 
sequences of length $n$ over $\z_v$ for which $\rho(b,t) = \kappa$. 
There exists a $m_{\rho},b_{\rho},c_{\rho} \in \R$ such that
  \[
    \sup_{\kappa} \left |\frac{\sqrt{2\pi b_{\rho}}r(\kappa)}{\binom{vl}{l,l,\ldots,l}}-c_{\rho} e^{(\kappa-m_{\rho})^2/b_{\rho}} \right | = o(1). 
  \]
\end{theorem}
\begin{proof}
Since the lattice $\Lambda = \z^v$ for either $Z = \{z\}$ or $Z = \{z \in \z_v^{t+2}:\; 
z(0) \neq x,\; z(t+1) \neq x,\; z(\iota) = x, 1 \leq \iota \leq t\}$, by 
Theorem~\ref{bender_thm}, in each case, there exists a mean $m_n \in \R^v$, 
covariance matrix $B_n \in \R^{v \times v}$ and constant $h \in \R$ such that
  \[
	  \sup_{k \in \N^v} \left | \frac{\sqrt{(2\pi)^v\det(B_n)}a_n(k)}{v^n} 
       - h e^{-(k-m_n)B_n(k-m_n)^T} \right | = o(1).
  \]
That is, the distribution of $a_n(k)$ is asymptotically normal. However this distribution 
is over all possible exponents $k$ which is equivalent to all circular sequences, 
not all balanced circular sequences as we desire. 

For $l \in \N$, let $k_{bal} = (l,l,\ldots,l,\kappa)$ and $a(\kappa) = a(k_{bal})$. 
The distribution of $a(\kappa)$ is a conditional distribution of the $a(k)$ and 
thus also satisfies a local limit theorem \cite[Proposition~3.13]{10.2307/20461449}. 
When $Z = \{z\}$, $a(\kappa) = t(\kappa)$ and when $Z = \{z \in \z_v^{t+2}:\; 
z(0) \neq x,\; z(t+1) \neq x,\; z(\iota) = x, 1 \leq \iota \leq t\}$, 
$a(\kappa) = r(\kappa)$, and the theorem is proved.
\end{proof}
We note that the conclusions of Theorem~\ref{thm_normal_random} hold when taken over sequences whose length is equivalent to $\beta \Mod v$ where $\beta \neq 0$ and where the numbers of any two symbols differ by at most 1.  For any fixed $\beta$, restricting $a_n(k)$ to the $k$ corresponding to these sequences is also a conditional distribution of the $a(k)$ and thus also satisfies a local limit theorem.  We are interested primarily in the case when $\beta = 0$ and in the next section we compute the mean and variance, whose existence is established in Theorem~\ref{thm_normal_random}, for this case.

\subsubsection{Mean and variance computation} \label{mean_variance}

In principle the mean and variance of this normal distribution could be computed 
by extracting the largest eigenvalue of $T$ and evaluating its derivatives. This 
is straightforward for any fixed values of $v$, $t$ and $Z$ but a closed form 
using this strategy is elusive. We compute the mean and variance directly, starting with tuples. 

We recall that for sequence $s$, we denote the length $(j-i)$ subsequence $(s(i), 
s(i+1), \ldots s(j-1))$ by $\subseq{s}{i}{j}$. This can be applied to any sequence 
including $\sigma$ and tuples $z$. We also recall that $|s|_a$ is the number of 
symbols $a$ that appear in $s$.  For a given $\sigma \in B(v,n)$ and $z \in \z_v^t$, let
\[
  I_{i,j} = \begin{cases} 1 & \mbox{if }  \sigma(i+_{n}\iota) 
= z(\iota) \mbox{ and } \sigma(j+_{n}\iota) = z(\iota), \; 0 \leq \iota < t; \\
  0 & \mbox{otherwise.} \end{cases}
\]
We observe that $I_{i,i}I_{j,j} = I_{i,j}$ and if $j < i+t$ and 
$z(\iota) \neq z(\iota+j-i)$ for any $0 \leq \iota < j-i$, then $I_{i,j} = 0$ for 
any $\sigma$. We will use a related indicator variable
\[
  \mathcal{I}_k = \begin{cases} 1 & \mbox{if }  z(\iota+k) = z(\iota),\; 0 \leq \iota < t-k; \\
  0 & \mbox{otherwise.} \end{cases}
  \]

In this section, we write $l = n/v$ to represent the number of times any symbol 
appears in a balanced sequence of length $n$ over $\z_v$. We also use $(n)_m$ 
to denote the descending product $ \prod_{i=0}^{m-1} (n-i) = n (n-1) \cdots (n-m+1)$. 

\begin{theorem}\label{E_lambda}
Let $z \in \z_v^t$ be a t-tuple over $\z_v$. The expected value and variance of $\lambda(z)$ 
taken over all balanced sequences of length $n$ over $\z_v$ are
  \begin{align*}
    E(\lambda(z)) =& \frac{\prod_{a=0}^{v-1} (l)_{|z|_a} }{(n-1)_{t-1}}, \\
    \VAR(\lambda(z)) =& \frac{\prod_{a=0}^{v-1} (l)_{|z|_a} }{(n-1)_{t-1}} - \left(\frac{\prod_{a=0}^{v-1} (l)_{|z|_a} }{(n-1)_{t-1}} \right)^2 + \frac{\prod_{a=0}^{v-1} (l)_{2|z|_a} }{(n-1)_{2t-2}} + % \\ &  
    2\sum_{k=1}^{t-1} \mathcal{I}_{k} \frac{\prod_{a=0}^{v-1} (l)_{|z|_a+|\subseq{z}{0}{k}|_a} }{(n-1)_{t+k-1}}.
\end{align*}
\end{theorem}

\begin{proof}
Let $\ell = n/v$. The number of balanced sequences of length $n$ over $\z_v$ is
  \[
    \binom{n}{\ell, \ell, \ldots, \ell}.
  \]
For any given $0 \leq i < n$, the number of sequences with 
$\sigma(i+_{n}\iota) = z(\iota)$ for all $0 \leq \iota < t$ is
  \[
    \binom{n-t}{\ell-|z|_0,\ell-|z|_1,\ldots,\ell-|z|_{v-1}}.
  \]
Therefore, we have 
\begin{align*}
E(\lambda(z)) &= E\left(\sum_{i=0}^{n-1} I_{i,i} \right)
               = \sum_{i=0}^{n-1} E(I_{i,i}) % \\ &
               = \frac{n\binom{n-t}{\ell-|z|_0,\ell-|z|_1,\ldots,
                 \ell-|z|_{v-1}}}{\binom{n}{\ell, \ell, \ldots, \ell}} 
               = \frac{\prod_{a=0}^{v-1} (l)_{|z|_a} }{(n-1)_{t-1}}.
\end{align*}

To compute the variance of $\lambda(z)$ we start by computing the second 
moment of $\lambda(z)$. Using properties of the indicator $I_{i,j}$ we obtain
\begin{align*}
E(\lambda(z)^2) &= E\left(\left(\sum_{i=0}^{n-1} I_{i,i}\right)
     \left(\sum_{j=0}^{n-1} I_{j,j}\right)\right) 
   = E\left(\left(\sum_{i=0}^{n-1} I_{i,i}^2 \right) + 
     \left(\sum_{i \neq j} I_{i,i}I_{j,j}\right)\right) \\
 & = E\left(\left(\sum_{i=0}^{n-1} I_{i,i} \right) + 
     \left(\sum_{i \neq j} I_{i,j}\right)\right) 
   = E(\lambda(z)) + \sum_{i \neq j} E(I_{i,j}). 
\end{align*}
The variance is then
\begin{align*}
  \VAR(\lambda(z)) &= E(\lambda(z)^2) - (E(\lambda(z)))^2\\
                   &= E(\lambda(z)) + \sum_{i \neq j} E(I_{i,j}) - (E(\lambda(z)))^2. 
\end{align*}
We examine the cross terms $\sum_{i \neq j} E(I_{i,j})$ more closely. We can assume without loss of generality that $0 \leq i < j < n$. If $j \geq i+t$ and $i+n \geq j+t$, then the two occurrences of $z$ at positions $i$ and $j$ do not overlap so $I_{i,j} = 1$. We calculate that the expected value is
\begin{align*}
 E(I_{i,j}) &= \frac{\binom{n-2t}{\ell-2|z|_0,\ell-2|z|_1,\ldots, 
               \ell-2|z|_{v-1}}}{\binom{n}{\ell,\ldots,\ell}} % \\ &
             = \frac{\prod_{a=0}^{v-1} (l)_{2|z|_a} }{(n)_{2t}}. 
\end{align*}

Next suppose that $j < i+t$ and $i+n \geq j+t$ or $j \geq i+t$ and $i+n < j+t$. 
If $j \geq i+t$ and $i+n < j+t$, we can rotate $\sigma$ and swap the roles of $i$ 
and $j$ to be in the first case. In these cases the copies of $z$ at positions 
$i$ and $j$ overlap. This overlap is possible, i.e. $I_{i,j} =1$,  when $z(\iota) = z(\iota+j-i)$ 
for $0 \leq \iota < t-(j-i)$ otherwise $I_{i,j}=0$. If overlap is possible, then
\begin{align*}
E(I_{i,j}) &= \frac{\binom{n-(t+j-i)}{\ell-|z|_0-|\subseq{z}{0}{j-i}|_0,\ell-|z|_1-|\subseq{z}{0}{j-i}|_1,
     \ldots, \ell-|z|_{v-1}-|\subseq{z}{0}{j-i}|_{v-1}}}{\binom{n}{\ell,\ldots,\ell}} % \\ &
   = \frac{\prod_{a=0}^{v-1} (l)_{|z|_a+|\subseq{z}{0}{j-i}|_a} }{(n)_{t+j-i}}. 
\end{align*}

Finally, if $j < i+t$ and $i+n < j+t$, then $n \leq 2t-2$ and $\sigma$ is 
completely determined by $z$, $i$ and $j$. Hence $I_{i,j} =1$ if $z(\iota) = z(\iota+j-i)$ for $0 \leq \iota < t-(j-i)$ and $z(\iota) = z(\iota+n+i-j)$ for $0 \leq \iota < t-(n+i-j)$ and $I_{i,j} = 0$ otherwise. Since $\sigma$ is balanced we must have that $l = |z|_a + |z|_{j-i}^{n+i-j}|_a$ for each $a$.  If overlap is possible in a balanced $\sigma$, then   
\[
E(I_{i,j}) = \frac{1}{\binom{n}{\ell,\ldots,\ell}} = \frac{(\ell!)^v}{n!}.
\]

For $n > 2t-2$ this case can never happen and we ignore this term.

% When $n > 2t-2$, 
Adding all these terms we get that
\begin{align*}
  \VAR(\lambda(z))  &= \frac{\prod_{a=0}^{v-1} (l)_{|z|_a} }{(n-1)_{t-1}} 
  - \left(\frac{\prod_{a=0}^{v-1} (l)_{|z|_a} }{(n-1)_{t-1}} \right)^2 % \\ &
  + n(n-2t+1) \frac{\prod_{a=0}^{v-1} (l)_{2|z|_a} }{(n)_{2t}} % \\ & 
  + 2n\sum_{k=1}^{t-1} I_{i,i+k} 
    \frac{\prod_{a=0}^{v-1} (l)_{|z|_a+|\subseq{z}{0}{k}|_a} }{(n)_{t+k}} \\
  &= \frac{\prod_{a=0}^{v-1} (l)_{|z|_a} }{(n-1)_{t-1}} - 
     \left(\frac{\prod_{a=0}^{v-1} (l)_{|z|_a} }{(n-1)_{t-1}} \right)^2 
     + \frac{\prod_{a=0}^{v-1} (l)_{2|z|_a} }{(n-1)_{2t-2}} + % \\ &
     2\sum_{k=1}^{t-1} \mathcal{I}_{k} 
     \frac{\prod_{a=0}^{v-1} (l)_{|z|_a+|\subseq{z}{0}{k}|_a} }{(n-1)_{t+k-1}}.
\end{align*}

\end{proof}

These exact expressions for the expectation and the variance have two problems: 
it is not easy to understand their magnitude even to first order simply by 
looking at the equations; they depend on the number of each symbol in $z$ and on 
the possibility of overlap.  To get a feeling for the magnitude of the expected 
value and variance of $\lambda$ in a way which does not depend on the 
particular $z$, we calculate bounds and use SageMath \cite{sagemath} to extract 
the leading terms. The following identities and approximations, whose proofs 
are reasonably straightforward, are used several times in our estimations.
\begin{lemma}\label{ident_approx}
\begin{align}
(x)_y &= x^y\left (1 - \frac{y(y-1)}{2x} + \frac{y(y-1)(y-2)(3y-1)}{24x^2} 
       - \frac{y^2(y - 1)^2(y - 2)(y - 3)}{48x^3} \right. \nonumber\\
      &\qquad\qquad\left .+ \frac{y(y-1)(y-2)(y-3)(y-4)(15y^3-30y^2+5y+2)}{5760x^4} 
       + O\left(\frac{1}{x^5} \right) \right), \label{x_y} \nonumber \\
\frac{1}{(x)_y} &= x^{-y} \left(1 + \frac{y(y-1)}{2x} + \frac{y(y-1)(y+1)(3y-2)}{24x^2} 
       + \frac{y^2(y-1)^2(y+1)(y+2)}{48x^3} \right. \nonumber\\
      &\qquad\qquad\left .+ \frac{y(y-1)(y+1)(y+2)(y+3)(15y^3-15y^2-10y+8)}{5760x^4} 
       +O\left(\frac{1}{x^5}\right) \right).  \nonumber
\end{align}
\end{lemma}
We also need to know how to maximize and minimize descending products with a 
fixed number of terms.
\begin{lemma}\label{min_max_x_y}
Let $s_i \geq 1$ for $0 \leq i < v$ and $\sum_{i=0}^{v-1} s_i = s$.  Then
$\prod_{i=0}^{v-1} (l)_{s_i}$ is minimized when one $s_i = s-v+1$ and all 
others are one, and is maximized when $|s_i -s_j| \leq 1$ for all $i$ and $j$.
\end{lemma}
\begin{proof}
Consider $a > b+1$, then
%   \begin{align*}
$$ (l)_a(l)_b = (l)_{a-1}(l-(a-1))(l)_b % \\
              < (l)_{a-1}(l-b)(l)_b % \\
              = (l)_{a-1}(l)_{b+1}. $$
%   \end{align*}
Hence, the product is minimized when all $s_i=1$ but one, that has to be 
$s-v+1$ to fulfill the conditions on $s$. The product is maximized when 
the $s_i$'s are as balanced as possible. 
\end{proof}
We conclude using Lemma \ref{min_max_x_y} that if $s = qv+r$, where 
$0 \leq r < v$, then
% \begin{equation}
$$ (l)_{s-v+1}l^{v-1} \leq \prod_{i=0}^{v-1} (l)_{s_i} 
   \leq ((l)_{q})^{v-r}((l)_{q+1})^{r}.$$ 
% \end{equation}

We are interested in computing the bounds with a final error term of $O(1/n)$, but to do that for the variance we need to compute one additional term in the expected value.  For a lower bound on the expectation, we compute
\[
  E(\lambda(z)) \geq   \frac{n}{v^t}\left( 1 + \frac{A_1(t,v)}{2n} + \frac{B_2(t,v)}{24n^2} \right ) +O\left(\frac{1}{n^2}\right), 
\]
where
\begin{align*}
  A_1(t,v) &= - (t^2 - 2tv + v^2 - t)(v - 1), \\
  B_1(t,v) &= (3(v - 1)t^4 - 2(6v^2 - v - 1)t^3 + 3(6v^3 + 2v^2 - v + 1)t^2 \\
           &- 2(6v^4 + 3v^3 - 3v^2 + v + 1)t + 3v^5 + v^4 - 2v^3)(v-1).
\end{align*}

To compute an upper bound we set $t = qv+r$ and use $0 \leq r\leq v-1$ to eliminate $r$.  We compute that
\[
  E(\lambda(z)) \leq \frac{n}{v^t} \left(1 + \frac{t(v-1)}{2n} +  \frac{(3(v+1)t^2 -2(v + 1)t + (3v^2 + 5v -6)v)(v - 1)}{24n^2}\right) +  O\left(\frac{1}{n^2}\right).
\]

We observe that although $E(\lambda(z))$ depends on the number of times each symbol 
appears in $z$, we determine upper and lower bounds on $E(\lambda(z))$ that are 
independent of $z$. We also note that if $t< v$, then using that $q=0$ and $r=t$, 
the upper bound can be improved to

\[
\frac{n}{v^t} \left(1 + \frac{t(t-1)}{2n} + \frac{(3t - 2)(t + 1)(t - 1)t}{24n^2} \right)
+ O\left(\frac{1}{n^2}\right).  
\]

In order to derive bounds on the variance, we first observe that the first two terms
\[
  E(\lambda(z)) - E(\lambda(z))^2 = E(\lambda(z))(1- E(\lambda(z)))
\] 
are bounded below by a lower bound on $E(\lambda(z))$ minus the product of an upper 
and lower bound. Similarly they are bounded above by an upper bound minus the same 
product. We only keep terms that allow a final error term of the order $O(1/n)$. 
Next we bound the terms that come from $\sum_{i \neq j} E(I_{i,j})$. In exactly 
the same manner as for the estimation of the expected value, we get

\[
  \frac{\prod_{a=0}^{v-1} (l)_{2|z|_a} }{(n-1)_{2t-2}} \geq  
  \frac{n(n-2t+1)}{v^{2t}}\left( 1 + \frac{A_2(t,v)}{2n} + 
  \frac{B_2(t,v)}{24n^2} \right ) +O\left(\frac{1}{n}\right),
\]
where
\begin{align*}
  A_2(t,v) &= - (4t^2 - 4tv + v^2 - 2t)(v - 1), \\ 
  B_2(t,v) &= (48(v-1)t^4 - 16(6v^2-v-1)t^3 + 12(6v^3 + 2v^2 - v + 1)t^2 \\ 
           &- 4(6v^4 + 3v^3 - 3v^2 + v + 1)t  + 3v^5+v^4-2v^3+2)(v-1).
\end{align*}
For the upper bound we get
\[
 \frac{\prod_{a=0}^{v-1} (l)_{2|z|_a} }{(n-1)_{2t-2}} \leq \frac{n(n-2t+1)}{v^{2t}} 
 \left(1 + \frac{2t(v-1)}{2n} +  \frac{(12(v+1)t^2 -4(v + 1)t 
 + (3v^2 + 5v + -6)v)(v - 1)}{24n^2}\right) +  O\left(\frac{1}{n}\right).
\]
We note that the term inside the parentheses is the same as that of the 
expectation but with $t$ replaced by $2t$.

For a lower bound on the fourth term of the variance, we take $I_{i,i+k}=0$ and thus this term 
is eliminated from our calculation of a lower bound of the variance. There are 
$z$ that cannot overlap themselves so $I_{i,i+k}=0$ does happen. For the upper 
bound we assume that $I_{i,i+k}=1$ and substitute $t+k$ for $t$ in the expression 
for the expectation and then sum the resulting polynomial divided by $v^{t+k}$. 
We use
\begin{align*}
  \sum_{k=1}^{t-1} \frac{1}{v^k} &= \frac{v^t-v}{(v-1)v^t}, \\
  \sum_{k=1}^{t-1} \frac{k}{v^k} &= \frac{v^{t+1}-(v^2-v)t - v}{(v-1)^2v^t}, \\
  \sum_{k=1}^{t-1} \frac{k^2}{v^k} &=   \frac{(v^{2} + v)v^{t} - (v-1)^2vt^2 -  2(v - 1)vt - (v^2 + v)}{(v-1)^3v^{t}},
\end{align*}
and obtain the upper bound
\begin{align*}
  2\sum_{k=1}^{t-1} I_{i,i+k} \frac{\prod_{a=0}^{v-1} (l)_{|z|_a+|\subseq{z}{0}{k}|_a} }{(n-1)_{t+k-1}} &= 2\sum_{k=1}^{t-1} I_{i,i+k} \frac{n}{v^{t+k}}\left (1 + \frac{(t+k)(v-1)}{2n} \right) + O\left(\frac{1}{n}\right)\\
&\leq \frac{n}{(v-1)v^{2t}} \left( 2(v^{t} - v) +   \frac{((v-1)t+v)v^t -2(v-1)vt  - v}{n}\right) + O\left(\frac{1}{n}\right).
\end{align*}

Putting all this together we get lower and upper bounds
\begin{align*}
  \VAR(\lambda(z)) &\geq \frac{n}{v^{2t}} \left( \frac{A_3(t,v)}{2} + \frac{B_3(t,v)}{24n} \right) + O\left(\frac{1}{n}\right), \\
    \VAR(\lambda(z)) &\leq \frac{n}{v^{2t}} \left( \frac{A_4(t,v)}{2} + \frac{B_4(t,v)}{24n} \right) + O\left(\frac{1}{n}\right),
\end{align*}
where
\begin{align*}
  A_3(t,v) &= (2v^t -3(v - 1)t^2 + 2(v^2 - v - 2)t + 2),\\
  A_4(t,v) &= \frac{(2(v + 1)v^t + (v-1)^2t^2 - 2(v-1)(v^2-v  +2)t 
            + v^4 - 2v^3 + v^2 - 2v - 2) }{(v-1)}, \\
% \end{align*}
% % {\small
% \begin{align*}
  B_3(t,v) &= -12(t^2 -2vt -t + v^2)(v-1)v^{t} + 45(v-1)^2t^4 -4(21v-26)(v+1)(v-1)t^3 \\
  &\quad + 6(9v^3+v^2-17v-14)(v-1)t^2 -12(v^4-2v^2-4v-2)(v-1)t -(3v-1)(v+6)(v-1)v,\\
  B_4(t,v) &= \left( 12(v^2t - t + 2v)v^t  - 3(v-1)^3t^4 + 4(v-1)^2(v+1)(3v-2)t^3  
  - 18(v^2-v+1)(v+2)(v-1)^2t^2  \right.\\
  &\quad \left. +12(v^3+v^2+2)(v^2-v-1)(v-1)t - (3v^6-5v^5-v^4+5v^3-2v^2+24)v \right)/(v-1).
\end{align*} % }

When $z$ cannot overlap with itself at any point the variance has a nicer form.
\begin{corollary}
  If $\subseq{z}{i}{t} \neq \subseq{z}{0}{t-i}$ for all $1 \leq i < t$, then 
  \[
 \VAR(\lambda(z)) =  \frac{\prod_{a=0}^{v-1} (l)_{|z|_a} }{(n-1)_{t-1}} - \left (  \frac{\prod_{a=0}^{v-1} (l)_{|z|_a} }{(n-1)_{t-1}} \right )^2 +   \frac{\prod_{a=0}^{v-1} (l)_{2|z|_a} }{(n-1)_{2t-2}}.
  \]
\end{corollary}

\begin{proof}
  If $\subseq{z}{i}{t} \neq \subseq{z}{0}{t-i}$ for all $1 \leq i < t$, then 
  $I_{i,j}=0$ for all $i+1 \leq j < i+t$ and $i-t < j \leq i-1$, with $I_{i,j} = 1$ 
  for all other $j$.
\end{proof}
This does not affect our lower bound since for the lower bound we assumed $I_{i,i+k}=0$ but when $z$ cannot overlap with itself the upper bound simplifies to
\[
    \VAR(\lambda(z)) \leq \frac{n}{2v^{2t}} \left(\frac{A_5(t,v)}{2}  + \frac{B_5(t,v)}{24n} \right ) +O\left(\frac{1}{n}\right),
\]
where
\begin{align*}
  A_5(t,v) & = 2v^t + (v - 1)t^2 - 2(v^2-v+2)t + (v^2-2v+2)(v+1), \\
  B_5(t,v) &= 12(v-1)tv^t - 3(v-1)^2t^4 + 4(v-1)(v+1)(3v-2)t^3 - 18(v^2-v+1)(v+2)(v-1)t^2 \\
  &\quad  +12(v^4+v^3-v^2+2)(v-1)t - (3v-2)(v+1)(v-1)v^3.
\end{align*}

Next we concentrate on runs. Over random balanced sequences, runs are easier 
to work with than tuples and we can provide more exact results.  Let $\rho(b,t)$ 
be the number of runs of $b$s of length $t$. Let $z_1 = (a_l,b,b,\ldots,b,a_r), 
z_2 = (c_l,b,\ldots,b,c_r) \in \z_v^{t+2}$ be two runs of symbol $b$ of 
length $t$, where $a_l,a_r,c_l,c_r$ are all different from $b$. Let
\[
  J_{i,z_1,j,z_2} = 
  \begin{cases} 
        1 & \mbox{if }  \sigma(i+_{n}\iota) = z_1(\iota) 
            \mbox{ and } \sigma(j+_{n}\iota) = z_2(\iota), \; 0 \leq \iota < t+2; \\  
        0 & \mbox{otherwise.} \end{cases}
\]
We are concerned mainly with the case that $n$ is large, so we assume that 
$n \geq 2vt$. Since $a_l,a_r,c_l,c_r \neq b$, we have that $J_{i,z_1,j,z_2} = 0$ 
if $i < j \leq i+t$ or $j < i \leq j+t$. If $i+t+1 < j < n+i-t-1$, then 
$J_{i,z_1,j,z_2} = 1$. When $j = i+t+1$, we have $J_{i,z_1,j,z_2} = 1$ if 
and only if $a_r = c_l$. Similarly, when $i = j+t+1$, we have 
$J_{i,z_1,j,z_2} = 1$ if and only if $a_l = c_r$.

\begin{theorem}
The expected number and variance of runs of length $t$ of a fixed symbol $b$ is
\begin{align*}
     E(\rho(b,t)) &= \frac{(l(v-1)-1)(v-1)l(l)_t}{(n-1)_{t+1}}, \\   
  \VAR(\rho(b,t)) &= \frac{(l(v-1)-1)(v-1)l(l)_t}{(n-1)_{t+1}} 
      - \left(\frac{(l(v-1)-1)(v-1)l(l)_t}{(n-1)_{t+1}}\right)^2 \\
    & + \frac{(v-1)l(l)_{2t}(l(v-1)-1)^2(l(v-1)-2)}{(n-1)_{2t+2}}.
\end{align*}
\end{theorem}
\begin{proof}
There are $v-1$ runs of length $t$ that start and end with the same symbol, 
and $(v-1)(v-2)$ that start and end with different symbols. Recalling that 
$l=n/v$ and using the same ideas as in Theorem \ref{E_lambda}, we get               
\begin{align*}
   E\left(\sum_{a_l,a_r \neq b}\sum_{i=0}^{n-1} J_{i,z_1,i,z_1}\right) 
&= (v-1)\frac{l(l-1)(l)_t}{(n-1)_{t+1}} + (v-1)(v-2) \frac{l^2(l)_{t}}{(n-1)_{t+1}}\\  
&= \frac{(l(v-1)-1)(v-1)l(l)_t}{(n-1)_{t+1}}. 
\end{align*}
For the variance we begin by using the fact that $J_{i,z_1,i,z_2}J_{j,z_1,j,z_2} 
= J_{i,z_1,j,z_2}$. Since $a_l,a_r,c_l,c_r \neq b$, runs can only overlap when 
$j = i+t+1$ or $j=n+i-t-1$. So $J_{i,z_1,j,z_2}=0$ whenever $n+i-t-1 < j < i+t+1$ 
and $j \neq i$. We obtain
\begin{align*}
  E\left(\left(\sum_{0 \leq i < n}\sum_{a_l,a_r \neq b} J_{i,z_1,i,z_1}\right)^2\right)
= & \sum_{0 \leq i < n \atop  j=i} 
% \sum{\begin{array}{c} 0 \leq i < n \\ j=i\end{array}}
    \sum_{a_l=c_l \neq b \atop a_r=c_r\neq b}
% \sum_{\begin{array}{l}a_l=c_l \neq b\\a_r=c_r\neq b\end{array}} 
    E(J_{i,z_1,i,z_1}) +  
    \sum_{0 \leq i < n \atop j= i+t+1}
%    \sum_{\begin{array}{c}0 \leq i < n \\ j= i+t+1\end{array}} 
    \sum_{a_l,a_r,c_l,c_r \neq b \atop a_r = c_l}
% \sum_{\begin{array}{c}a_l,a_r,c_l,c_r \neq b \\ a_r = c_l\end{array}} 
    E(J_{i,z_1,j,z_2}) \\
+ & \sum_{0 \leq i < n \atop j= n+i-t-1}
% \sum_{\begin{array}{c}0 \leq i < n \\ j= n+i-t-1\end{array}} 
    \sum_{a_l,a_r,c_l,c_r \neq b \atop a_l = c_r}    
% \sum_{\begin{array}{c}a_l,a_r,c_l,c_r \neq b \\ a_l = c_r\end{array}} 
    E(J_{i,z_1,j,z_2}) + 
    \sum_{{0 \leq i < n \atop j\geq i+t+2} \atop j\leq n+i-t-2 }
% \sum_{\begin{array}{c}0 \leq i < n \\ j\geq i+t+2 \\ j\leq n+i-t-2 \end{array}}
    \sum_{a_l,a_r,c_l,c_r \neq b} E(J_{i,z_1,j,z_2}).
\end{align*}

The first term is just $E(\rho(b,t))$. The second and third terms are 
equal, with exactly three non-$b$ values to be chosen, $\{a_l,a_r,c_r\}$ 
and $\{a_l,a_r,c_l\}$, respectively. These three values can all be the same, 
two be the same and one differ or all three be distinct. In the case that 
there are two different symbols, there are three ways to choose the pair 
of positions which get the same symbol. Thus, the sum of the second and 
third terms is  
\begin{align*}
  & \quad \sum_{0 \leq i < n \atop j= i+t+1}
    \sum_{a_l,a_r,c_l,c_r \neq b \atop a_r = c_l}
    E(J_{i,z_1,j,z_2}) +  \sum_{0 \leq i < n \atop j= n+i-t-1}
    \sum_{a_l,a_r,c_l,c_r \neq b \atop a_l = c_r}    
    E(J_{i,z_1,j,z_2}) \\
 =& \quad \frac{2n(l)_{2t}}{(n)_{2t+3}} \left((v-1) (l)_3 
    + 3(v-1)(v-2)l^2(l-1) + (v-1)(v-2)(v-3)l^3 \right) \\
 =& \quad \frac{2(v-1)l(l)_{2t}(l(v-1)-1)(l(v-1)-2)}{(n-1)_{2t+2}}.
\end{align*}

For the last term we must consider the number of different symbols in and partitions 
of $\{a_l,a_r,c_l,c_r\}$, obtaining
\begin{align*}
  & \quad  \sum_{{0 \leq i < n \atop j\geq i+t+2} \atop j\leq n+i-t-2 }
    \sum_{a_l,a_r,c_l,c_r \neq b} E(J_{i,z_1,j,z_2})\\
 =& \quad \frac{n(n-2t-3)(v-1)(l)_{2t}}{(n)_{2t+4}} \\
  & \quad \left((l)_4 + 4(v-2)l(l)_3 + 3(v-2)(l)_2(l)_2
    +6(v-2)(v-3)l^2(l)_2+(v-2)(v-3)(v-4) l^4\right) \\
 =& \quad \frac{(v-1)l(l)_{2t}(l(v-1)-1)(l(v-1)-2)(l(v-1)-3)}{(n-1)_{2t+2}}.
\end{align*}

Putting all pieces together, we obtain that the variance is
\begin{align*}
 & \VAR(\rho(b,t)) \\
% =& \frac{(l(v-1)-1)(v-1)l(l)_t}{(n-1)_{t+1}} 
%  - \left(\frac{(l(v-1)-1)(v-1)l(l)_t}{(n-1)_{t+1}}\right)^2 \\
% +& \frac{n(v-1)(l)_{2t}}{(n)_{2t+3}}\left((l)_3(l-1) + (v-2)l^2(l-1)(7l-5) 
%    + 2(v-2)(v-3)l^3(3l-2) + (v-2)(v-3)(v-4) l^4 \right) \\
=& \frac{(l(v-1)-1)(v-1)l(l)_t}{(n-1)_{t+1}} 
 - \left(\frac{(l(v-1)-1)(v-1)l(l)_t}{(n-1)_{t+1}}\right)^2 
 + \frac{(v-1)l(l)_{2t}(l(v-1)-1)^2(l(v-1)-2)}{(n-1)_{2t+2}}.
\end{align*}
\end{proof}
We observe that these expressions do {\em not} depend on the symbol $b$ 
in the run. To understand how these expressions behave, we compute the 
leading terms using SageMath \cite{sagemath}. For the expected value we get

\[
E(\rho(b,t)) = \frac{n(v-1)}{v^{t+2}} \left((v-1) 
    - \frac{(v-1)^2t^{2} -(v+3)(v-1)t +2}{2n}\right) + O\left (\frac{1}{n}\right).
\]
% \begin{multline}
% E(\rho(b,t)) = \frac{n(v-1)}{v^{t+2}} \left((v-1) 
%     - \frac{(v-1)^2t^{2} -(v+3)(v-1)t +2}{2n}  \right . \\
% + \left . \frac{3(v-1)^3t^4 - 2(v-1)^2(5v+11)t^3 + 3(v-1)(3v^2+6v+19)t^2 - 2(v^3 - v^2 - 7v + 31)t - 24 }{24n^2} \right)+ O\left (\frac{1}{n^2}\right)
% \end{multline}
For the variance, we get
\[
\VAR(\rho(b,t)) = \frac{n(v-1)}{v^{2t+4}} \left( A_6(t,v) 
     + \frac{B_6(t,v)}{2n}\right) + O\left(\frac{1}{n}\right)   
\]
where
\begin{align*}
  A_6(t,v) =& \; (v-1)\left(v^{t+2} - (v-1)^3t^2 + 2(v-1)^2t - (v-1)(v+1)\right),\\
  B_6(t,v) =& \; v^{t+2}(-(v-1)^2t^2+(v-1)(v+3)t-2) + 3(v-1)^5t^4 -4(v-1)^4(v+4)t^3 \\
    & + (v-1)^3(v^2+8v+31)t^2 - 2(v-1)^2(v^2+2v+13) + 8(v-1).
 \end{align*}

\section{Sequences from ElGamal}\label{sec_elgamal}

At this point we have determined the balance, periodicity, tuple balance and run properties of $\z_v$ sequences reduced from random permutations of $\zpx$ by the remainder operator.  In this section we consider permutations of $\zpx$ produced from the ElGamal 
function and compare their properties to those from random permutations. Just as in Section~\ref{sec_random}, let $p$ be a prime. Our primary interest is in the case when $v \mid p-1$ but we state some of our results more generally. If $g \in \zpx$ is primitive, that is $g$ is a generator of the 
multiplicative group $\zpx$, then $\gamma = (g^{i-1}\rem p)_{i=1}^{p-1}$ 
is a permutation of $\zpx$ which we call a {\em $\zp$-ElGamal sequence} or 
just {\em ElGamal sequence}. We denote $\gamma_v = ((g^{i-1}\rem p)\rem v)_{i=1}^{p-1}$ 
and call it an {\em ElGamal sequence modulo $v$}. The permutation $\gamma$ 
and sequence $\gamma_v$ are properly functions of $p$, $v$ and $g$ but we 
omit explicitly using functional notation because it risks confusion with 
our other notations. We recall that we fix the underlying set of $\z_p$ to 
be $\{0, 1, \ldots , p-1\}$, that of $\z_{p}^*$ to be $\{1, 2,\ldots, p-1\}$ 
and determine $x \rem v \in \z_v$ by considering $x \in \z_p^*$ to be the 
corresponding element in $\{0, 1, \ldots, p-1\} \subset \z$.

%%%%%%%%%%%%%%%%%%%%%%%%%%%%%%%%%%%%%%%%%%%%%%%%%%%%%%%%%%%%%%%%%%%%%%%%%%%%% 
\subsection{Balance property and period lengths of ElGamal sequences modulo $v$}\label{sec_balance_elgamal}

Since $\gamma$ is a permutation of $\zpx$, Proposition \ref{balancedproperty} 
shows that $\gamma_v$ is always balanced and exactly balanced if and only if $v \mid p-1$. 

We have proved that, if $\pi \in \zpx$ is a permutation picked at random, then 
the probability that $\pi_v$ has maximal period length approaches 1 very quickly 
as $p$ grows. Here we prove that ElGamal sequences have maximal period always.

\begin{theorem}\label{thm_period_elgamal}
For $p$ an odd prime and $g$ primitive in $\zp$, let 
$\gamma_v = ((g^{i-1}\rem p)\rem v)_{i=1}^{p-1}$. The ElGamal sequence 
$\gamma_v$ has maximal period length $N = p-1$.
\end{theorem}
\begin{proof}
If $p \equiv \alpha \Mod v$ with $\alpha \neq 1$, then Lemma~\ref{lemma_near_balanced_period} shows that $\gamma_v$ has period $p-1$. If $\alpha =1$, suppose that $\gamma_v$ has period $N < p-1$.  Thus we have that
  \[
    g^{i+N}\rem p \equiv_v g^i \rem p 
  \]
for all $0 \leq i < N$. Let $g' = g^N\rem p $. Setting $i=0$ shows that 
$g' \equiv_v 1$.  Let $p = kg'+r$ with $0 \leq r < g'$, set $x = k+1$ and let 
$j = \log_g(x)$. Setting $i =j$ gives that $x \equiv_v xg'\rem p$. Since 
$p < xg' < 2p$, $xg'\rem p = xg'-p$.  Thus
\[ xg'\rem  p \equiv_v xg'-p \equiv_v xg' - 1.\]
Now $x \equiv_v xg' \rem  p$ gives $x \equiv_v xg' -1$, or equivalently 
$x (g'-1) \equiv_v 1$. This contradicts that $g' \equiv_v 1$. Thus,  
$\gamma_v$ has period $p-1$.  
\end{proof}

%%%%%%%%%%%%%%%%%%%%%%%%%%%%%%%%%%%%%%%%%%%%%%%%%%%%%%%%%%%%%%%%%%%%%%%%%%%%% 
\subsection{$t$-tuple balance and run properties of ElGamal sequences modulo $v$}
For $z \in \z_v^t$ a $t$-tuple over $\z_v$ and 
$\gamma_v = ((g^{i-1}\rem p)\rem v)_{i=1}^{p-1}$ the ElGamal sequence with $p$ 
prime and $g$ primitive, let
\[
  \lambda(z) = \cardinality{\{i \in [1,p-1]: \gamma_v(i+_{n} \iota) = z(\iota),\; 0 \leq \iota < t\}},
\]
where we recall that $i+_{n} \iota = (i+\iota) \rem n$. 
This counts the number of times a fixed $t$-tuple appears in $\gamma_v$. 

For $b \in \z_v$, $t \in \mathbb{N}$, let
\[
  \rho(b,t) = \cardinality{\{i \in [1,p-1]: \gamma_v(i-_{n} 1) \neq b,\; \gamma_v(i+_n t) \neq b,\;  
  \gamma_v(i+_n\iota) = b,\; 0 \leq \iota < t\}}.
\]
This counts the number of runs of length $t$ of a fixed symbol.

\subsubsection{$t$-tuples}\label{subsec_tuples}

\begin{theorem}\label{elgamal_tuple}
  Let $p$ be an odd prime, $g$ be primitive in $\zp$ and $\gamma_v$ the corresponding ElGamal sequence modulo $v$.  If
  \[
    p = q g^{t-1} + r
  \]
  with $0 \leq r < g^{t-1}$, then
  \[
     \left \lfloor \frac{g}{v} \right \rfloor^{t-1} 
     \left \lfloor\frac{q}{v} \right \rfloor 
\leq \lambda(z) 
\leq \left \lceil \frac{g}{v} \right \rceil^{t-1} 
     \left(\left \lfloor\frac{q}{v} \right\rfloor +1 \right).
  \]
\end{theorem}

\begin{proof}
Let $X = \{x \in [1,p-1]: (g^i x) \rem p \equiv_v z(i) \mbox{ for all }
0 \leq i < t \}$. Because the sequence $\gamma$ is a permutation of $[1,p-1]$,
  \[
    \lambda(z) = \cardinality{X}.
  \]
We partition this set into parts whose sizes we can bound fairly precisely. 
Let $p \equiv \alpha \Mod v$.  Since $p$ is prime, $\alpha$ is invertible in $\z_v$.  For $0 \leq i < t$ define $c_i = g^{i}z(0) - z(i)$. Let
  \[
    D = \{d \in \z^t:d_0 = 0,  d_i \equiv_v \alpha^{-1}c_i\mbox{ and } 
          gd_{i-1} \leq d_i < g(d_{i-1}+1) \mbox{ for } 0 < i < t  \}.
  \]
  For $d \in D$, let\[
    X_{d} = \left \{x \in \z: x \equiv_v z(0),\; \frac{d_{i}p}{g^i} 
    \leq x < \frac{(d_{i}+1)p}{g^i}, \mbox{ for } 0 \leq i < t \right \}.
  \]
  
  We show that 
  \[
    X = \bigcup_{d \in D} X_d.
  \]
Suppose that $x \in X$ and thus $(g^i x) \rem p \equiv_v z(i)$ for all 
$0 \leq i < t$. Let us define $q_i$ and $r_i$ by $g^ix = q_ip + r_i$ 
with $0 \leq r_i < p$. Remark that $q_0 = 0$ and that 
$r_i = g^ix - q_ip = (g^ix) \rem p \equiv_v z(i)$ for all $0 \leq i < t$. 
First we observe that
  \[
    \frac{q_{i}p}{g^i} \leq x < \frac{(q_{i}+1)p}{g^i}  
  \]
for all $0 \leq i < t$ and $x \equiv_v (g^0x)\rem p \equiv_v z(0)$. 
Thus, if $\bar{q}=(q_0,\ldots,q_{t-1})\in \z^t$, then % $x \in X_{(q_i)_0^{t-1}}$. 
$x \in X_{\bar{q}}$. We also have
  \begin{align*}
    q_i &\equiv_v \alpha^{-1}q_ip  
          = \alpha^{-1}(g^ix - r_i )
          \equiv_v \alpha^{-1}(g^iz(0) - z(i)) = \alpha^{-1}c_i.
  \end{align*}
Then,
  \begin{align*}
    q_i &= \frac{g^ix-r_i}{p} = \frac{g(g^{i-1}x)-r_i}{p}
          = \frac{g(q_{i-1}p+r_{i-1}) - r_i}{p} \\
        &= gq_{i-1} + g\frac{r_{i-1}}{p} -\frac{r_i}{p} 
          < g(q_{i-1}+1).
  \end{align*}
 We have
  \[
  gq_{i-1} = \frac{gq_{i-1}p}{p} \leq \frac{g(q_{i-1}p+r_{i-1})}{p} 
           = \frac{g(g^{i-1}x)}{p} = \frac{g^ix}{p} = q_i + \frac{r_i}{p}.
  \]
Since $gq_{i-1}$ and $q_{i}$ are both integers and $0\leq r_{i}/p < 1$ we 
have that $q_{i} \geq g q_{i-1}$. Thus
  \[
    g q_{i-1} \leq q_i < g(q_{i-1}+1).
  \]
Therefore $\bar{q}=(q_0,\ldots,q_{t-1}) \in D$.

We now show that for any $d \in D$, $X_{d} \subset X$. For $d \in D$, 
if $x \in X_{d}$, then
  \[
    d_{i}p \leq g^ix < (d_{i}+1)p 
  \]
for all $0 \leq i < t$ and $x \equiv_v z(0)$. Then we have
  \begin{align*}
    g^ix \rem  p &= g^ix-d_{i}p 
                \equiv_v g^ix-\alpha d_i 
                \equiv_v g^iz(0) - c_i 
                \equiv_v g^iz(0) - (g^iz(0) - z(i)) = z(i),
  \end{align*}
which shows that $x \in X$.

Thus 
  \begin{align*} 
    X = \bigcup_{d \in D}X_d % \\ % 
      = \bigcup_{d \in D} \left (\{x \equiv_v z(0) \} \bigcap 
        \left( \bigcap_{0 \leq i < t} \left \{\frac{d_ip}{g^i} 
        \leq x < \frac{(d_i+1)p}{g^i} \right\} \right) \right). \\     
  \end{align*}
However since $gd_{i-1} \leq d_i < g(d_{i-1}+1)$, the sets of the 
inner-most intersection are nested, thus
  \begin{align*}
    X = \bigcup_{d \in D}
    \left (   \{x \equiv_v z(0) \} \cap \left \{\frac{d_{t-1}p}{g^{t-1}} 
    \leq x < \frac{(d_{t-1}+1)p}{g^{t-1}} \right \}\right). 
  \end{align*}

The index set of the outer union is at least $\lfloor g/v \rfloor^{t-1}$ 
and at most $\lceil g/v \rceil^{t-1}$. We know that the real valued 
interval $[d_{t-1}p/g^{t-1},(d_{t-1}+1)p/g^{t-1})$ on the right hand 
side of the intersection has length $q + r/g^{t-1}$ and since $p$ is 
prime, $r > 0$.  The fewest number of integers in such an interval is 
$q$ and the largest number is $q+1$.  Thus the size of the intersection 
is at least $\lfloor q/v \rfloor$ and at most $\lceil (q+1)/v \rceil$. Finally there can be no integer strictly between $q/v$, and $(q+1)/v$ 
and thus $\lceil (q+1)/v \rceil = \lfloor q/v \rfloor + 1$.
\end{proof}

One of our motivations in studying ElGamal sequences modulo $v$, besides comparing them to random balanced sequences, is 
understanding how they behave with respect to randomness measures. If 
$g = mv$, then $\lfloor g/v \rfloor = \lceil g/v \rceil = m$ and the 
bounds on $\lambda(z)$ differ by at most $m^t$. So if $g$ is a multiple 
of $v$ and $m^t$ is relatively small, $\gamma_v$ is close to having the 
$t$-tuple balance property.  When $g = v$, then $m^t=1$ and the sequences 
have the $t$-tuple balance property. In this case we have
\[
  \left \lfloor\frac{q}{v} \right \rfloor \leq \lambda(z) \leq   
  \left \lfloor\frac{q}{v} \right\rfloor +1.
\]
In this case suppose that there are $n_{l}$ and $n_{u}$ tuples where 
$\lambda(z) = \left \lfloor q/v \right \rfloor,\; 
\left \lfloor q/v \right \rfloor +1$, respectively. The two equations
\begin{align*}
  n_{l} + n_{u} &=  v^t, \\
  n_{l}\left \lfloor \frac{q}{v} \right \rfloor + n_{u} \left(\left \lfloor \frac{q}{v} \right \rfloor +1 \right ) &= p-1,
\end{align*}
determine $n_{l}$ and $n_{u}$ exactly. We note that $n_l >0$; supposing 
otherwise, for a contradiction, that $n_{l} = 0$ gives 
$n_{u} = v^t$ and that $p = k v^t + 1$ for some $k$. Thus, $q = kv$ and 
$\lfloor q/v \rfloor = k$.   Plugging this back into the two linear equations gives $n_{u} = 0$ and the contradiction $n_l = v^t$. From the fact that $n_l > 0$, we can conclude that when $g=v$ there always exist tuples 
$z \in \z_v^t$ with $\lambda(z) = \lfloor q/v \rfloor$. The same 
considerations show that $n_u = 0$ if and only if $p \equiv 1 \pmod {v^t}$.

Artin's conjecture on primitive roots proposes that every integer $v$ 
which is not $-1$ nor a perfect square, is a primitive element for 
infinitely many primes. Hooley proved this conjecture under the 
generalized Riemann Hypothesis, however no unconditional proof exists 
\cite{MR966133,MR3011564}. It has been proved that there do exist 
$v$ for which this is true but no single specific value is known. 
For example, there are at most two primes which are primitive for 
only finitely many $p$. This implies that at least one of $2$, $3$, 
$5$ is primitive for infinitely many $p$. Another example is that 
there are at most three square-free numbers which are primitive for 
only finitely many $p$. For any such prime, setting $g=v$ gives an ElGamal sequence which satisfies the $t$-tuple balance property.

In general the relative errors introduced by the floors and ceilings are reduced when the numerators are large with respect to $v$, the denominator.  So when $g$ is not a multiple of $v$ the relative errors introduced by the floors and ceilings of $g/v$ are minimized when $g$ is as large as possible. However this only goes so far because as $g$ increases, $q$ decreases and the relative errors from $q/v$ will increase.  Rough calculation shows that when $g$ and $q$ are not  multiples of $v$, the bounds are tightest when $g \approx \sqrt[t]{p}$.

We conclude our study of $t$-tuples in ElGamal sequences giving a series
of related results and corollaries.

\begin{corollary}
Let  $p$ be an odd prime, $g$ primitive and $\gamma_v$ the corresponding 
ElGamal sequence modulo $v$.  If $p \geq v g^{t-1}$ and $g \geq v$, 
then $\lambda(z) > 0$ for all $z \in \z_v^t$.
\end{corollary}
\begin{proof}
If $g \geq v$, then $\lfloor g/v \rfloor \geq 1$. If $p \geq vg^{t-1}$, 
then $v \leq q$ and $\lfloor q/v \rfloor > 1$. Thus 
Theorem~\ref{elgamal_tuple} gives that $\lambda(z) \geq 1$ for any $z$.
\end{proof}
\begin{corollary}
Let $p$ be an odd prime, $g$ primitive and $\gamma_v$ the corresponding 
ElGamal sequence modulo $v$. If $\lambda(z) > 0$ for all $z \in \z_v^t$, 
then $g \geq v$ and $p \geq v^t+1$.
\end{corollary}
\begin{proof}
If $g < v$, then there are tuples $(c_i)_{i=0}^{t-1}$ for which there 
can not be any $d$ with $d_i \equiv_v \alpha^{-1}c_i$. In other words $D$ 
is empty. In fact there are $g^t$ tuples $(c_i)_{i=0}^{t-1}$ which admit 
a $d$ and $v^t-g^t$ tuples $(c_i)_{i=0}^{t-1}$ which do not. Thus it is 
necessary that $g \geq v$.

If $\gamma_v$ contains every $t$-tuple at least once, then the length $p-1$ of 
$\gamma_v$ must be at least $v^t$.  
\end{proof}
When $g = v$ these necessary and sufficient conditions coincide.
\begin{corollary}

Let  $p$ be an odd prime, $g=v$ primitive and $\gamma_v$ the corresponding 
ElGamal sequence modulo $v$.  Then $\lambda(z) > 0$ for all $z \in \z_v^t$ 
if and only if $p \geq v^t+1$
\end{corollary}

When $g < v$ we can precisely give which tuples $z \in \z_v^t$ cannot occur.
\begin{proposition}
Let $p$ be an odd prime, $g<v$ primitive and $\gamma_v$ the corresponding 
ElGamal sequence modulo $v$. Then for all $i$, $\gamma_v(i+1) \equiv_v 
g\gamma_v(i)- s$ for some $ 0 \leq s < g$.
\end{proposition}
\begin{proof}
Let $g^{i} = q p + r v + \gamma_v(i)$ with $0 \leq rv + \gamma_v(i) < p$.  
Also let $grv+g\gamma_v(i) = s p + r'v + \gamma_v(i+1)$ with $0 \leq r'v + \gamma_v(i+1) < p$. Since $rv + \gamma_v(i) < p$, we conclude that 
$0 \leq s < g$. Then, since $p \equiv_v 1$, we have
\begin{align*}
    \gamma_v(i+1) &= (g^{i+1} \rem p )\rem v % \\ &
                   = ( g(g^i) \rem  p ) \rem v % \\ &
                   = (g(qp+rv+\gamma_v(i)) \rem p) \rem v \\
                  &= (grv+g\gamma_v(i)) \rem  p) \rem  v % \\ &
                   = (grv+g\gamma_v(i) - sp )\rem v % \\ &
                   \equiv_v g\gamma_v(i) -s.
\end{align*}

\end{proof}

\subsubsection{Runs}
To count the number of runs we use a slight variation on Theorem~\ref{elgamal_tuple}.
\begin{theorem}\label{elgamal_run}
Let  $p$ be an odd prime, $g$ primitive and $\gamma_v$ the corresponding ElGamal sequence 
modulo $v$.  For $z \in \z_v^t$, let 
  \[
    \mu(z) = \cardinality{\{ i \in [1,p-1]:\; g^{i+j}\rem p \equiv_v z(j),\; 0 \leq j < t-1,
                    \; g^{i+t-1}\rem p \not\equiv_v z(t-1)\}}.
  \]
  If
  \[
    p = q g^{t-1} + r
  \]
  with $0 \leq r < g^{t-1}$, then
  \[
    \left \lfloor \frac{g}{v} \right \rfloor^{t-2} \left \lfloor \frac{(v-1)g}{v} \right \rfloor \left \lfloor\frac{q}{v} \right \rfloor \leq \mu(z) \leq  \left \lceil \frac{g}{v} \right \rceil^{t-2} \left \lceil \frac{(v-1)g}{v} \right \rceil\left(\left \lfloor\frac{q}{v} \right\rfloor +1 \right).
  \]
\end{theorem}
\begin{proof}
  Let $X = \{x \in [1,p-1]:\; (g^i x) \rem p \equiv_v z(i) \mbox{ for all } 0 \leq i < t-1,\; 
  (g^{t-1} x) \rem p \not \equiv_v z(t-1)\}$. Then $\mu(z) = \cardinality{X}$.
As in Theorem~\ref{elgamal_tuple}, we partition $X$. Let $p \equiv \alpha \Mod v$. For $0 \leq i < t$, define 
$c_i = g^{i}z(0) - z(i)$. Let
  \[
    D = \{d \in \z^t:d_0 = 0,  d_i \equiv_v \alpha^{-1}c_i,\; 1 \leq i < t-2,\; 
        d_{t-1} \not\equiv \alpha^{-1}c_{t-1}\mbox{ and } gd_{i-1} \leq d_i < g(d_{i-1}+1) \mbox{ for } 0 < i < t \}.
  \]
Just as before
  \begin{align*} 
    X &= \bigcup_{d \in D} \left( \{x \equiv_v z(0) \} \cap \left \{\frac{d_{t-1}p}{g^{t-1}} \leq x < \frac{(d_{t-1}+1)p}{g^{t-1}} \right\}\right).
  \end{align*}
The index set of the outer union is at least $\lfloor g/v \rfloor^{t-2}\lfloor (v-1)g/v \rfloor$ 
and at most $\lceil g/v \rceil^{t-2}\lceil (v-1)g/v \rceil$. The size of the intersection is at 
least $\lfloor q/v \rfloor$ and at most $\lfloor q/v \rfloor + 1$. Now, the remaining of the proof 
is as in Theorem~\ref{elgamal_tuple}.
\end{proof}

We can use this theorem to count runs in two different ways.  For $b \in \z_v$, let $\rho(b,t)$ be the number of runs of symbol $b$ of length $t$ in the sequence $\gamma_v$.  That is,
\[
  \rho(b,t)  = \cardinality{ \{ i \in [1,p-1]:\; g^i\rem p \not\equiv_v b, g^{i+t+1}\rem p \not\equiv_v b,\; 
  g^{i+j}\rem p \equiv_v b,\; 1 \leq j \leq t\}}.
\]
\begin{corollary}\label{elgamal_runs1}
  Let  $p$ be an odd prime, $g$ primitive and $\gamma_v$ the corresponding ElGamal sequence modulo $v$.  If
  \[
    p = q_t g^{t} + r_t,\quad p = q_{t+1} g^{t+1} + r_{t+1}
  \]
  with $0 \leq r_i < g^{i}$ for $i = t, t+1$, then
\begin{eqnarray*}
 & & \left \lfloor \frac{g}{v} \right \rfloor^{t-1}\left \lfloor \frac{(v-1)g}{v} \right \rfloor \left \lfloor \frac{q_t}{v}\right \rfloor - \left \lceil \frac{g}{v} \right \rceil^{t}\left \lceil \frac{(v-1)g}{v} \right \rceil \left \lceil \frac{q_{t+1}+1}{v} \right \rceil  \\
 & \leq & \rho(b,t) \leq \left \lceil \frac{g}{v} \right \rceil^{t-1}\left \lceil \frac{(v-1)g}{v} \right \rceil \left \lceil \frac{q_{t}+1}{v} \right \rceil - \left \lfloor \frac{g}{v} \right \rfloor^{t}\left \lfloor \frac{(v-1)g}{v} \right \rfloor \left  \lfloor \frac{q_{t+1}}{v} \right \rfloor.
\end{eqnarray*}
\end{corollary}
\begin{proof}
  The value of $\mu(b^{t+1})$ is the number of runs of $b$ of length at least $t$.  Thus $\rho(b,t) = \mu(b^{t+1})-\mu(b^{t+2})$ and the bounds follow.
\end{proof}

\begin{corollary}\label{elgamal_runs2}
  Let  $p$ be an odd prime, $g$ primitive and $\gamma_v$ the corresponding ElGamal sequence modulo $v$.  If
  \[
    p = q g^{t+1} + r
  \]
  with $0 \leq r < g^{t+1}$, then
  \[
    (v-1) \left \lfloor \frac{g}{v} \right \rfloor^{t} \left \lfloor \frac{(v-1)g}{v} \right \rfloor\left \lfloor\frac{q}{v} \right \rfloor \leq \rho(b,t) \leq  (v-1)\left \lceil \frac{g}{v} \right \rceil^{t}\left \lceil \frac{(v-1)g}{v} \right \rceil \left \lceil\frac{q+1}{v} \right\rceil.
  \]
\end{corollary}
\begin{proof}
  Since $\rho(b,t) = \sum_{a \neq b} \mu(ab^{t+1})$, the bounds follow. 
\end{proof}

Remarkably, either of these corollaries can give better bounds. For example if $p=1759$, 
$v=2$ and $t=3$, when $g=6$ the lower bound in Corollary~\ref{elgamal_runs1} is better, but for 
all other primitive $g$ in $\z_{1759}$, the lower bound in Corollary~\ref{elgamal_runs2} is better. 
If $p=1097$, $v=2$ and $t=3$, the upper bound from Corollary~\ref{elgamal_runs2} is better when $g=3,5$, the two bounds are equal if $g=6$ and the upper bound from Corollary~\ref{elgamal_runs1} is better for all other primitive $g$.  In experiments with $100 \leq p < 10000$, $v < 7$, $t < 6$ and a sampling of small values of $g$ (so the lower bounds are likely to be larger than 0), the lower bound from 
Corollary~\ref{elgamal_runs1} is better for $0-5\% $ of parameters sets $(p,g,v,t)$, Corollary~\ref{elgamal_runs2} is better for $75-85\%$ of parameters sets and the two lower bounds are equal for $16-19\%$ of parameter sets, with some variation for different $t$.  For the same range of parameters the upper bound from Corollary~\ref{elgamal_runs1} is better for $73-83\% $ of parameters sets $(p,g,v,t)$, Corollary~\ref{elgamal_runs2} is better for $0-8\%$ of parameters sets and the two upper bounds are equal for $15-20\%$ of parameter sets, again varying with $t$.

If $g$ is invertible in $\z_v$, then we can improve the counting for the runs.
\begin{theorem} \label{g_invertible}
Let $p$ be an odd prime, $g$ primitive and $\gamma_v$ the corresponding ElGamal sequence 
modulo $v$. If
  \[
    p = q g^{t+1} + r
  \]
  with $0 \leq r < g^{t+1}$, then
  \[
    \left \lfloor \frac{g}{v} \right \rfloor^{t-1} \left \lfloor \frac{(v-1)g}{v} \right \rfloor^2 \left \lfloor\frac{q}{v} \right \rfloor \leq \rho(b,t) \leq  \left \lceil \frac{g}{v} \right \rceil^{t-1}\left \lceil \frac{(v-1)g}{v} \right \rceil^2 \left \lceil\frac{q+1}{v} \right\rceil.
  \]
\end{theorem}
\begin{proof}
  The proof proceeds as in Theorem~\ref{elgamal_run} but with
\begin{eqnarray*}
    D &=& \{d \in \z^{t+2}:d_0 = 0,\;  d_1 \not \equiv_v \alpha^{-1}(g-1)b,\; 
            d_i \equiv_v \alpha^{-1}(g^{i-1}d_1 + (g^{i-1}-1)b),\; 2 \leq i \leq t,\\
      & & \; d_{t+1} \not\equiv \alpha^{-1}(g^td_1 + (g^t-1)b), \mbox{ and } 
             gd_{i-1} \leq d_i < g(d_{i-1}+1) \mbox{ for } 0 < i < t+2  \},
\end{eqnarray*}
  and
  \[
    X = \{x \in [1,p-1]:\; x \equiv_v b+(\alpha d_1-(g-1)b)g^{-1},\; d_ip \leq g^ix < (d_i+1)p,\; 0 \leq i < t+2\}.
  \]
\end{proof}
This is an improvement since we are not multiplying the errors introduced by 
$\lfloor g/v \rfloor$ and $\lceil g/v \rceil$ by $v-1$.

Certainly in very general terms the bounds given in Corollaries~\ref{elgamal_runs1} 
and~\ref{elgamal_runs2} behave as the run property asks: there are $v$ times more 
runs of length $i$ than runs of length $i+1$.  Just as in the case of $tuples$, 
when $g$ is divisible by $v$ the ElGamal sequence more closely matches the run 
property. When $g=v$, using Corollary~\ref{elgamal_runs2}, the total number of runs 
of length $t$ is between
\[
v(v-1)^2 \left \lfloor\frac{q_{t+1}}{v} \right \rfloor
\]
and
\[
v(v-1)^2 \left \lceil\frac{q_{t+1}+1}{v} \right\rceil,
\]
which differ by $v(v-1)^2$. When $g$ is not a multiple of $v$, the same pattern 
of relative accuracy as $g$ grows is observed for the bounds on runs as it was 
observed for the bounds on tuples.

Golomb's run postulate refers to the number of runs of a given length, for all possible symbols, that is, $\rho(t) = \sum_{b \in \z_v} \rho(b,t)$.  We have not found better bounds for $\rho(t)$ other than multiplying our bounds by $v$ although we expect that summing over all symbols improves the distribution about the mean by canceling out lower order terms of opposite sign.  To first order, all of our bounds give $\rho(b,t) = p/v^t$, giving $\rho(t) = p/v^{t-1}$.  Golomb's run postulate asks that $\rho(t) \approx v \rho(t+1)$ which our bounds show is true in first order terms.  

\section{Experimental results} \label{experiments}

We now present some experimental data on the behaviour of $\lambda(z)$, $\rho(b,t)$ 
and $\rho(t)$ in ElGamal sequences modulo $v$ with respect to our bounds and 
desired properties. For tuples and runs, we denote the lower and upper bounds 
as $lb$ and $ub$, respectively. 
For each $2 \leq v \leq 8$ we chose $100$ distinct pairs $(p,v)$ such that
$p > 1,000,000$ and $v | p-1$. Using the smallest 10 generators, we construct
the ElGamal sequences and directly calculate the number of occurrences of tuples 
of length $2 \leq t \leq 7$ for each pair $(p,v)$. We repeat the experiment
considering pairs $(p,v)$ with the additional condition where $v = g$.

\begin{figure}
  \includegraphics[width=8cm]{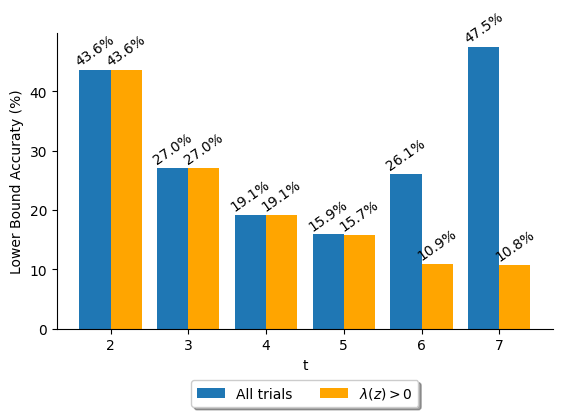} \includegraphics[width=8cm]{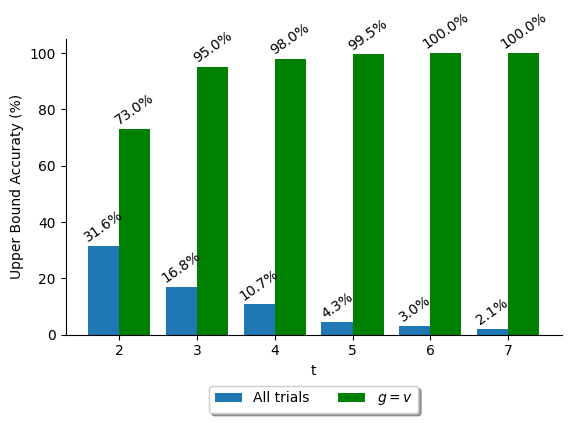}
  \caption{Percentage of trials with $z \in \z_v^t$ such that $\lambda(z)$ matches the lower bounds (left figure) and upper bounds (right figure). \label{bound_accuracy}}
\end{figure}

\begin{figure}
  \begin{center}
  \begin{tabular}{cc}
 \includegraphics[width=4cm]{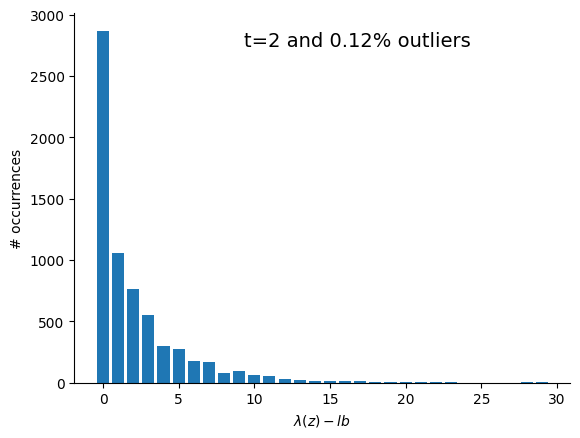} & \includegraphics[width=4cm]{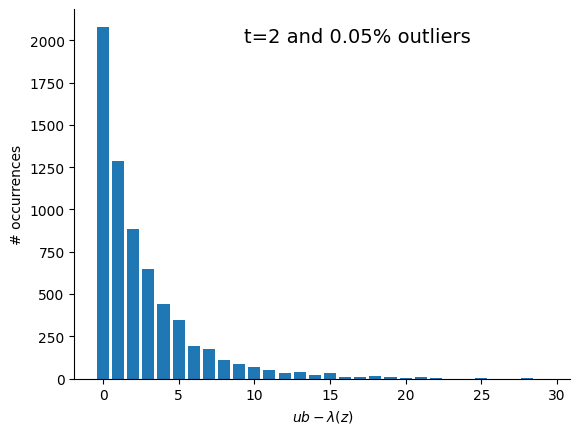} \\
 \includegraphics[width=4cm]{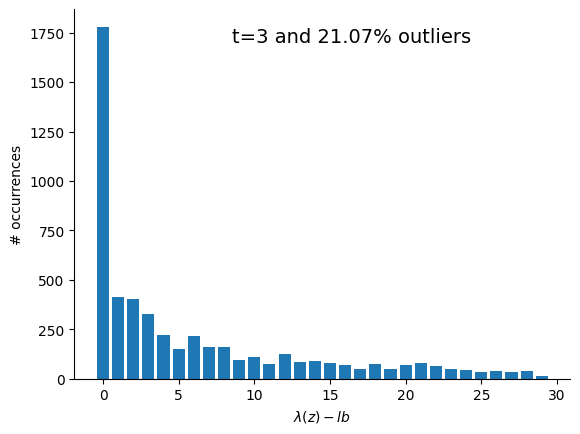} & \includegraphics[width=4cm]{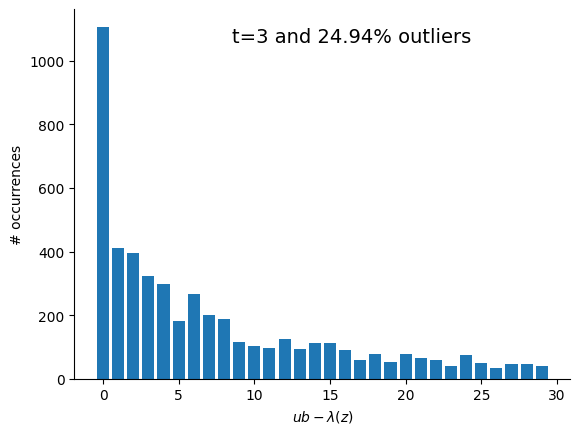} \\
 \includegraphics[width=4cm]{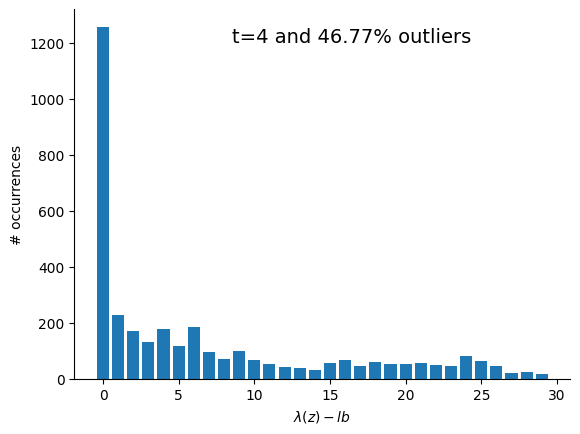} & \includegraphics[width=4cm]{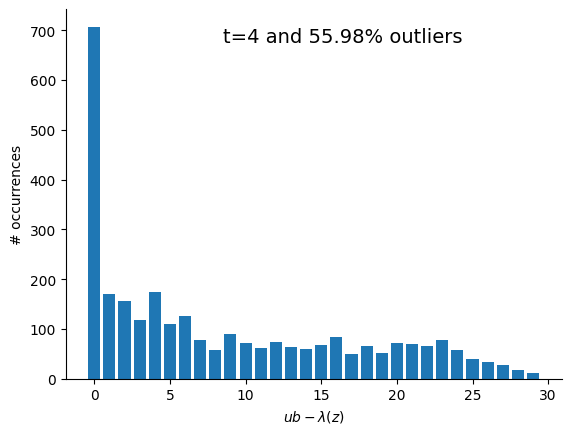} \\
 \includegraphics[width=4cm]{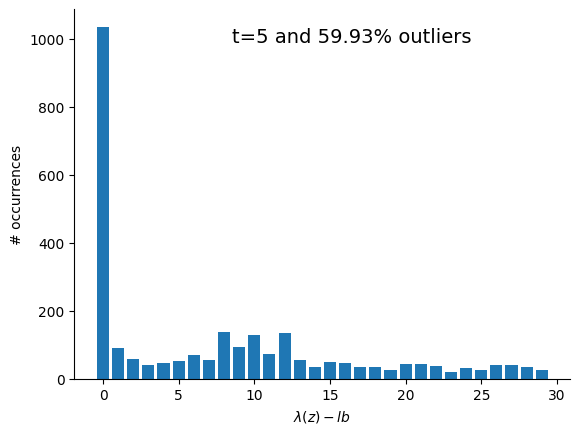} & \includegraphics[width=4cm]{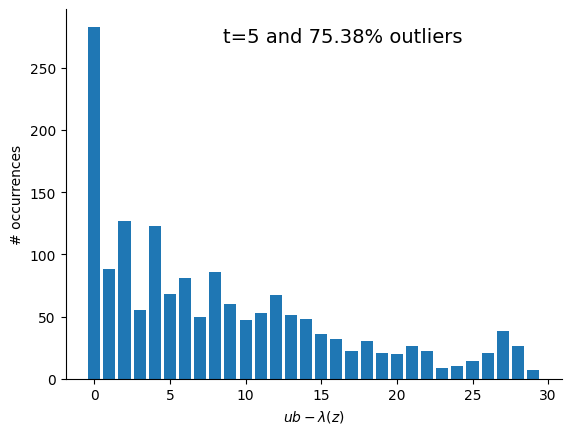} \\
 \includegraphics[width=4cm]{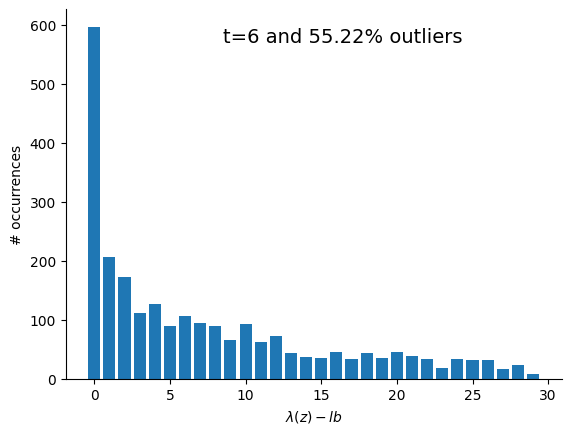} & \includegraphics[width=4cm]{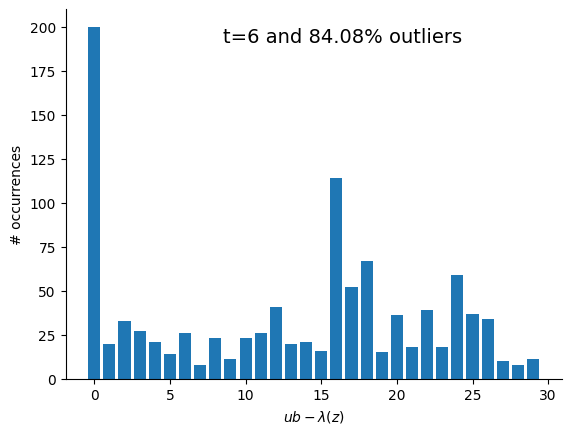} \\
    \includegraphics[width=4cm]{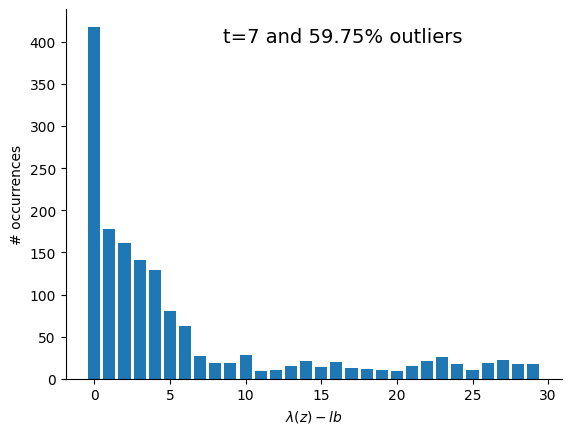} & \includegraphics[width=4cm]{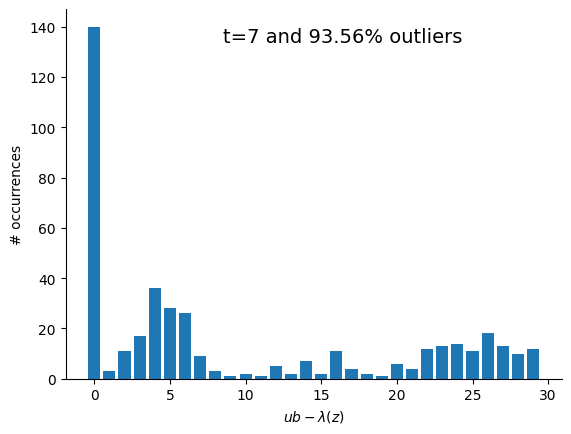}
  \end{tabular}
  \end{center}
  \caption{Distribution of gaps between extreme $\lambda(z)$ and the lower (left figures) and upper (right figures) bounds. \label{bound_dist}}
\end{figure}

Figure~\ref{bound_accuracy} shows the percentage of these trials where there 
exist $z \in \z_v^t$ such that $\lambda(z)$ matches the lower bound. The blue 
bars show the numbers over all the trials. For lower bounds the orange bars 
show only those trials when $\lambda(z) > 0 $ for all $z \in \z_v^t$; these 
are the cases we are most interested in because no tuples are missing. The 
upper bounds are never zero so this data is only shown for lower bounds. For 
upper bounds the green bars show the data for trials when $g=v$. In 
Section~\ref{subsec_tuples} we prove that when $g=v$, there always exist 
$z$ with $\lambda(z)$ matching the lower bound but the same is not true with 
respect to the upper bound so we report the data for upper bounds only. The 
green data show the percentage of the time $n_u>0$. Figure~\ref{bound_dist} 
shows the distribution of $\min\{\lambda(z):z \in \z_v^t\}-lb$ and 
$ub - \max\{\lambda(z):z \in \z_v^t\}$ when these values are small. The 
percentage of $z \in \z_v^t$ not appearing in the data is shown as the outliers 
in each case. For lower bounds the values ignore the trials where 
$\min\{\lambda(z)\} = 0$. These data show that our bounds are met with equality 
fairly frequently so we should expect any improvements to the bounds to contain 
conditions on the relationships between parameters, or depend on $z \in \z_v^t$. 
We can also see that even when the bounds are not tight, they are frequently 
very close.

\begin{figure}
  \includegraphics[width=8cm]{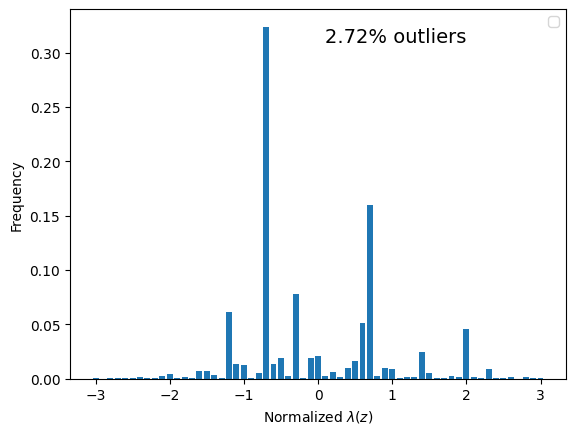} \includegraphics[width=8cm]{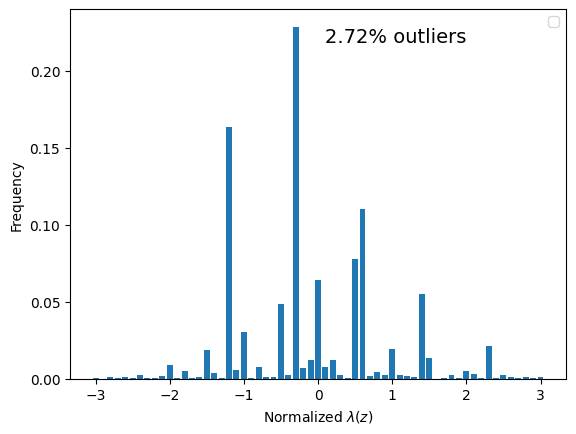}
  \caption{Distribution of normalized values of $\lambda(z)$ across all trials (left) and when $g=v$ (right). \label{normalized_lambda}}
\end{figure}

The main barrier to make a meaningful comparison of the behaviour of values 
of $\lambda(z)$ that occur in ElGamal sequences to $\lambda(z)$ occurring in 
sequences derived from random permutations (as described in 
Theorems~\ref{thm_normal_random} and~\ref{E_lambda}) is the sparseness of the 
data for any fixed parameters.  To address this we plot transformed values of $\lambda(z)$ over all experiments. 
The $\lambda(z)$ values are translated by $(p-1)/v^t$ so the mean of all values in 0.  Furthermore they are scaled by $\sqrt{v^t/(p-1)}$ which, from Theorem~\ref{E_lambda}, is the scaling factor required to make the distributions that exist from Theorem~\ref{thm_normal_random} coincide together with a variance of 1.  These data are shown in Figure~\ref{normalized_lambda} 
for all trials and for those trials where $g=v$. As can be seen the distributions
are not normal, but are not totally inconsistent with a constrained sampling from a normal distribution.  We think that for any definitive conclusion to be drawn about this behaviour of 
ElGamal sequences modulo $v$ when compared to sequences derived from random 
permutations, a much larger data set needs to be plotted.

\begin{figure}
  \includegraphics[width=8cm]{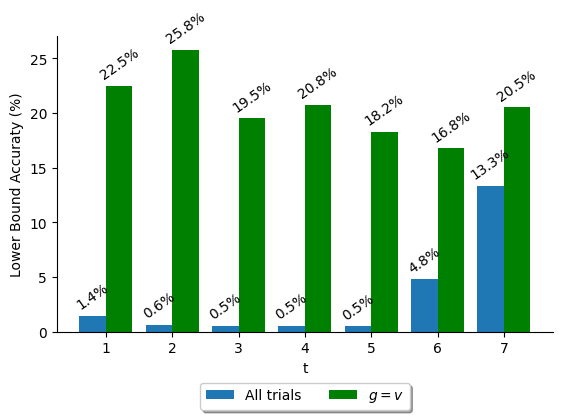} \includegraphics[width=8cm]{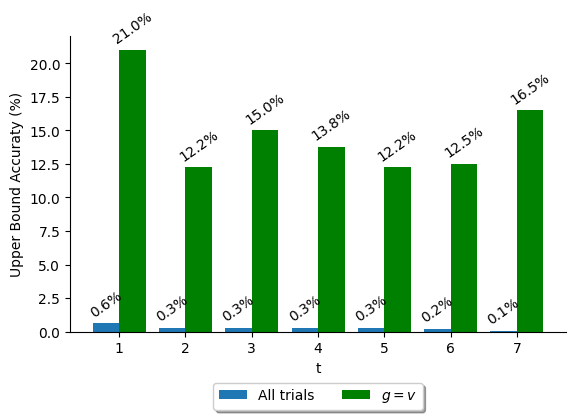}
  \caption{Percentage of trials with $b \in \z_v$ such that $\rho(b,t)$ matches the lower bounds of Corollary~\ref{elgamal_runs2} (left figure) and upper bounds of Corollary~\ref{elgamal_runs1} (right figure). \label{run_bound_accuracy}}
\end{figure}

\begin{figure}
  \begin{center}
  \begin{tabular}{cccc}
    \includegraphics[width=4cm]{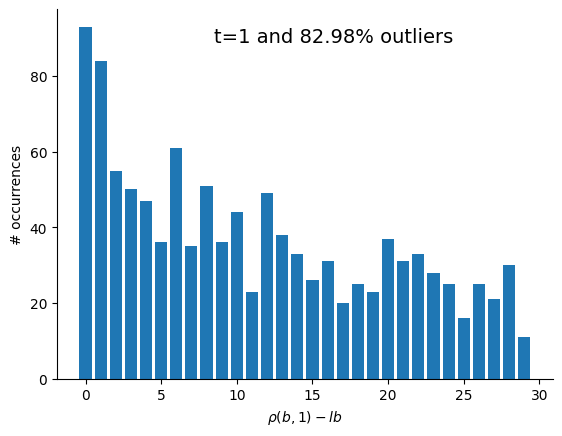} & \includegraphics[width=4cm]{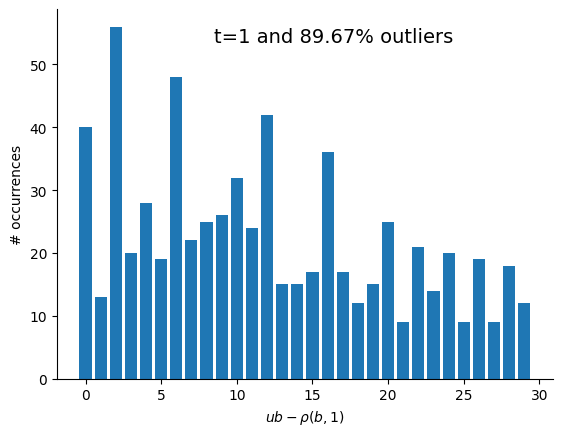} & \includegraphics[width=4cm]{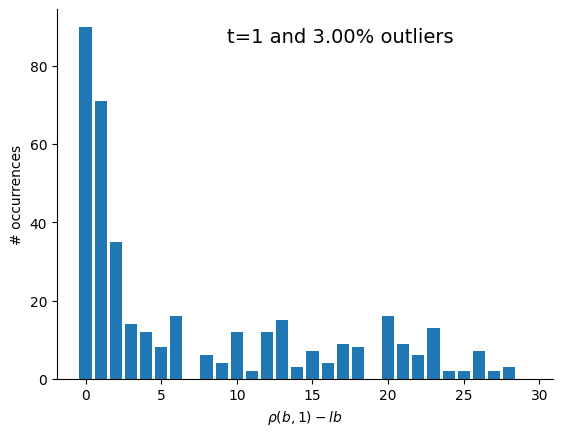} & \includegraphics[width=4cm]{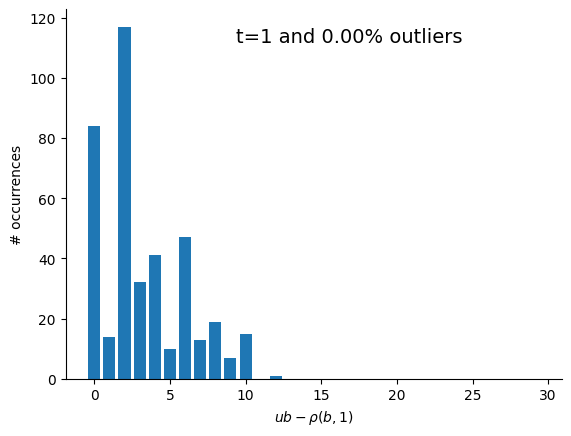} \\
    \includegraphics[width=4cm]{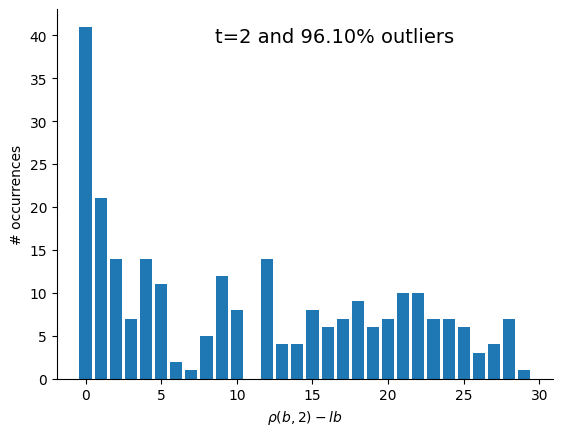} & \includegraphics[width=4cm]{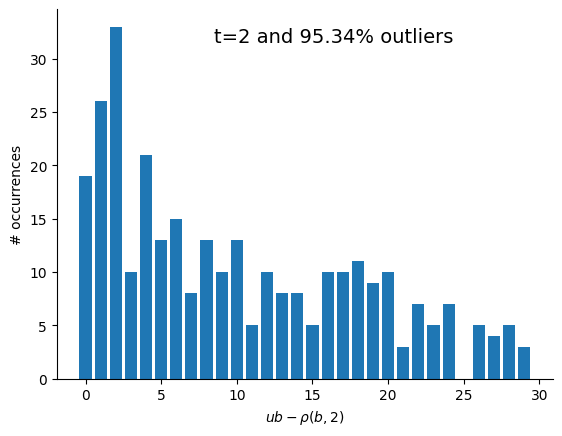} & \includegraphics[width=4cm]{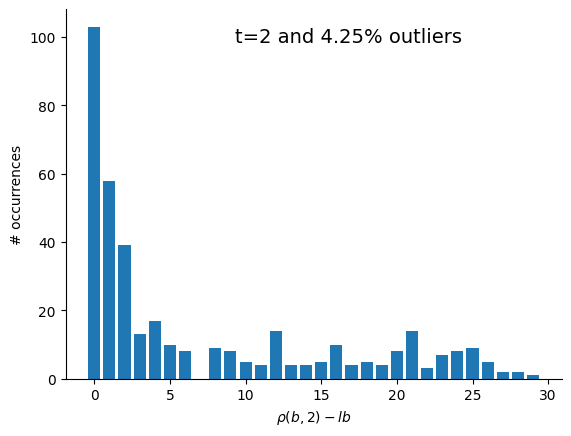} & \includegraphics[width=4cm]{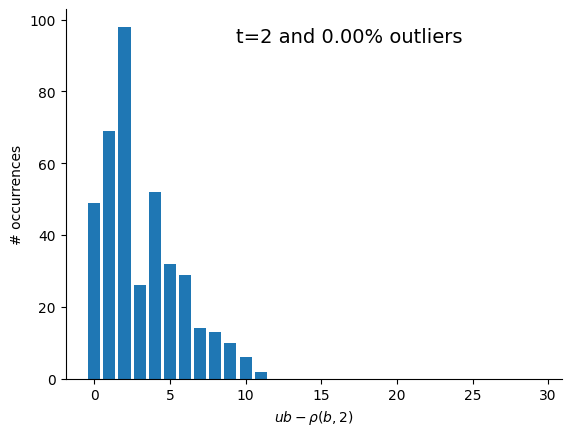} \\
    \includegraphics[width=4cm]{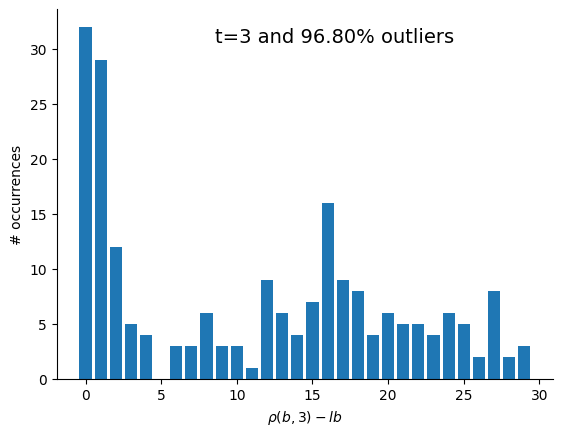} & \includegraphics[width=4cm]{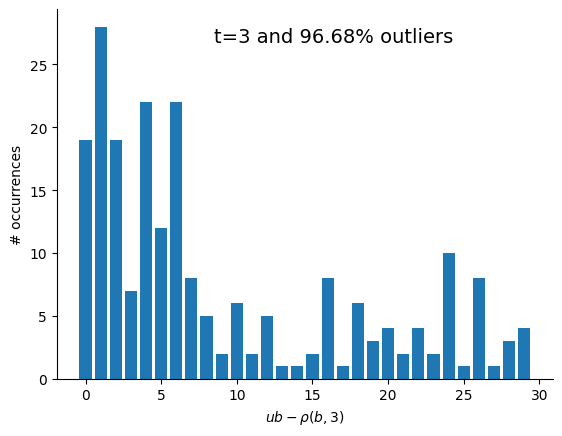} & \includegraphics[width=4cm]{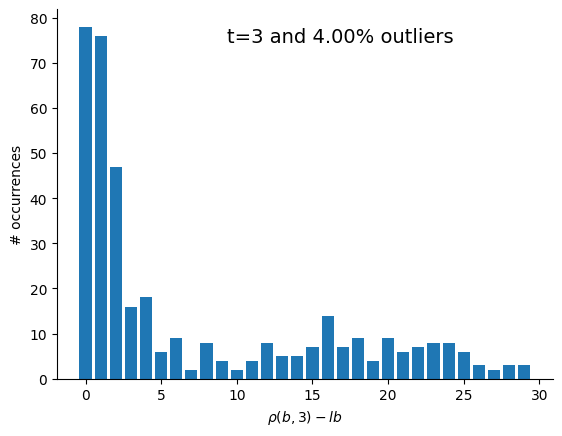} & \includegraphics[width=4cm]{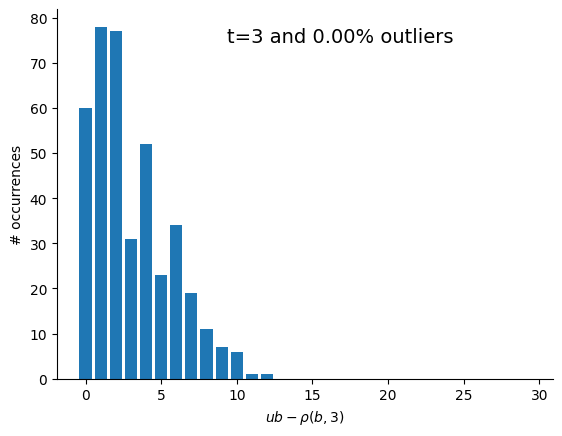} \\
    \includegraphics[width=4cm]{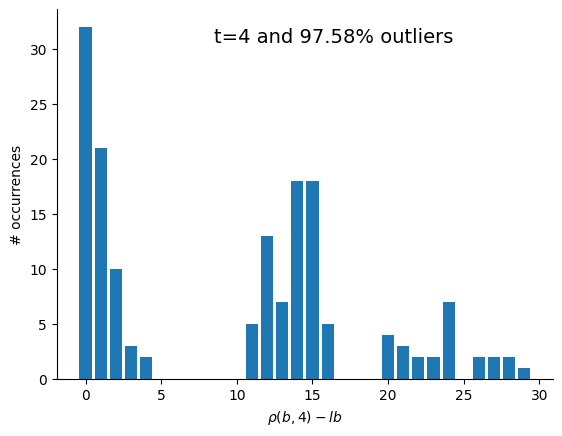} & \includegraphics[width=4cm]{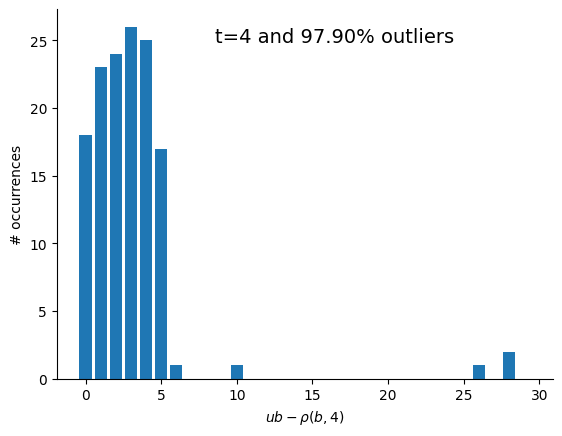} & \includegraphics[width=4cm]{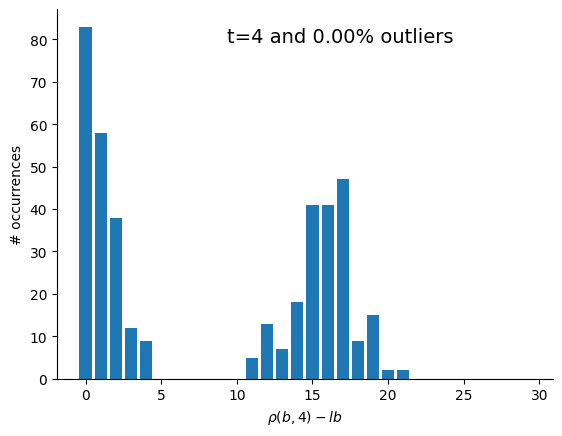} & \includegraphics[width=4cm]{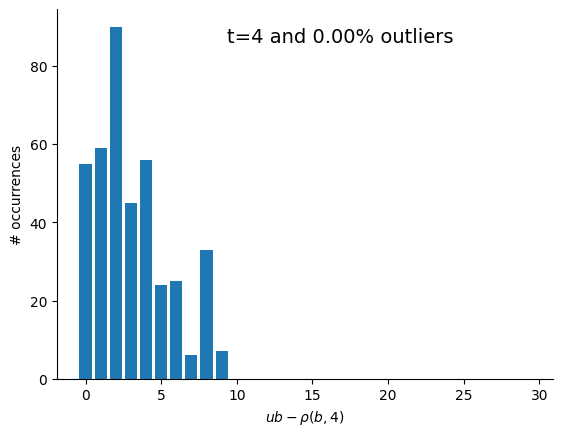} \\
    \includegraphics[width=4cm]{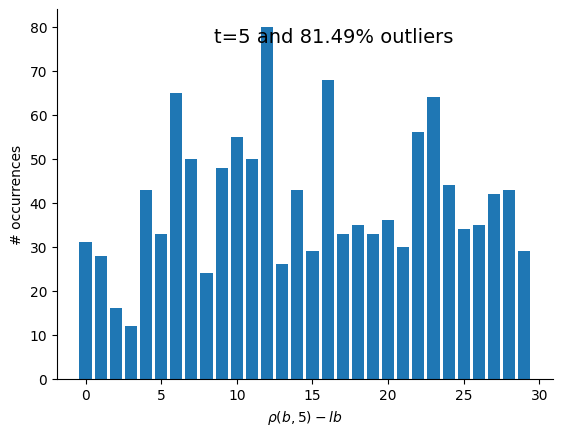} & \includegraphics[width=4cm]{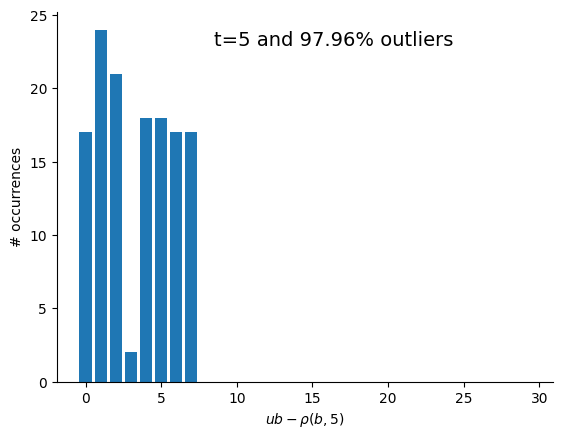} & \includegraphics[width=4cm]{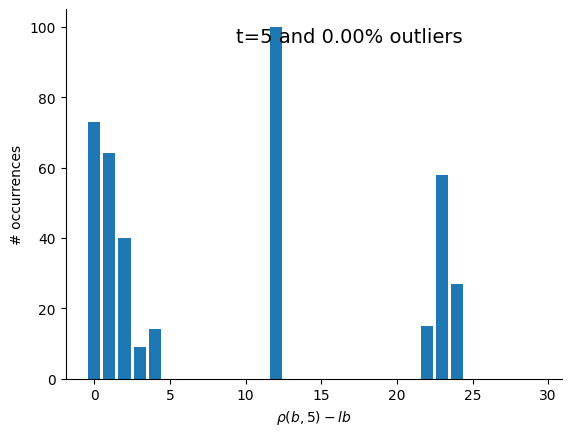} & \includegraphics[width=4cm]{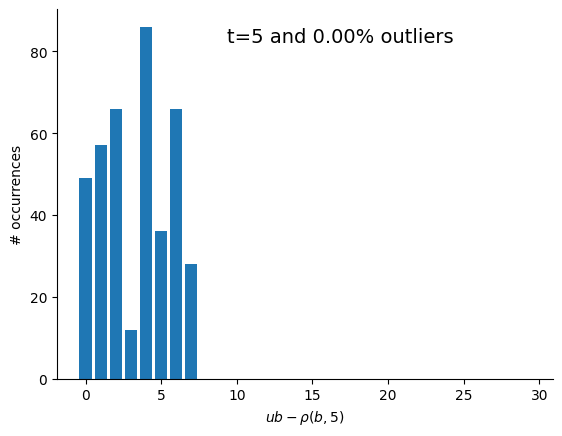} \\
    \includegraphics[width=4cm]{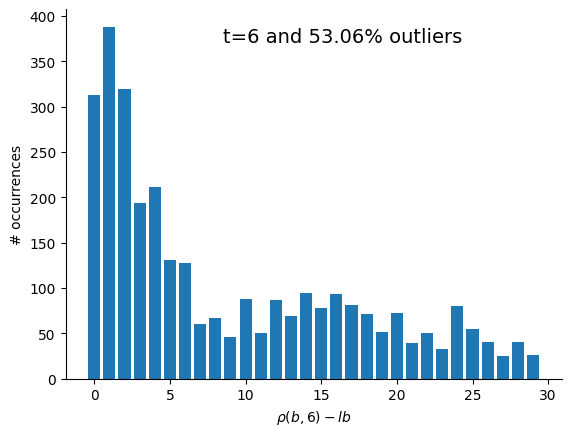} & \includegraphics[width=4cm]{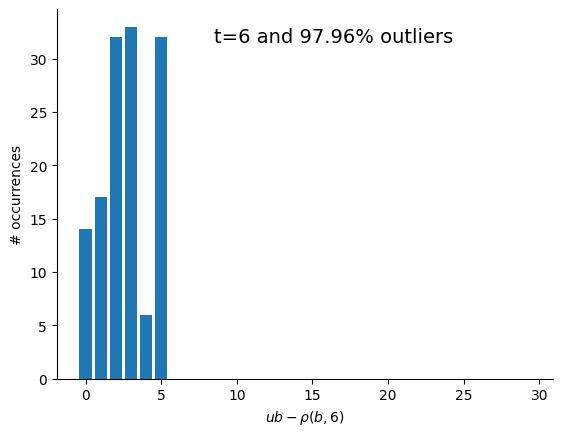} & \includegraphics[width=4cm]{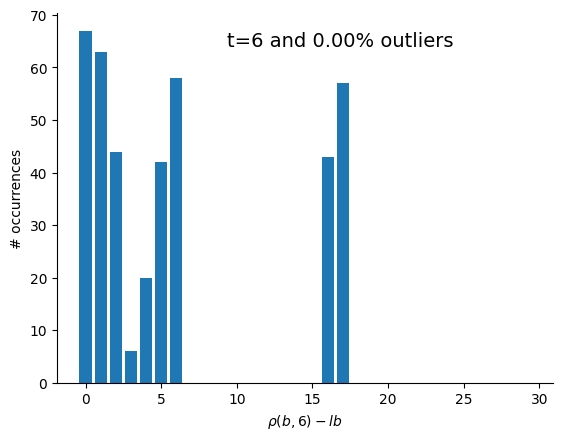} & \includegraphics[width=4cm]{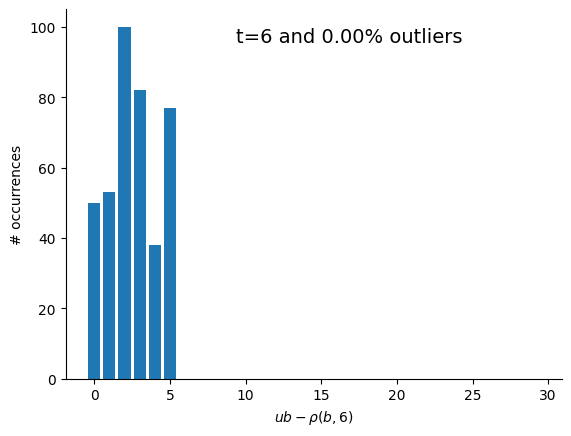} \\
    \includegraphics[width=4cm]{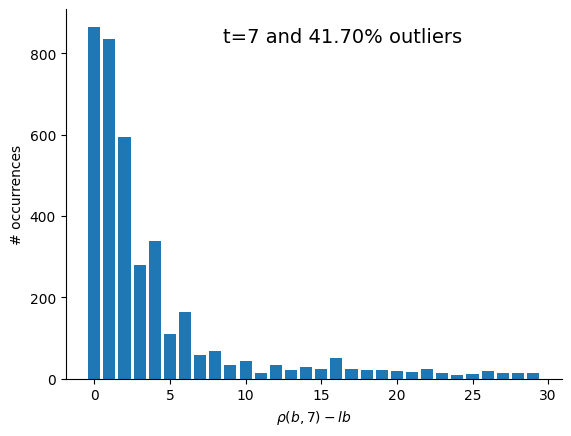} & \includegraphics[width=4cm]{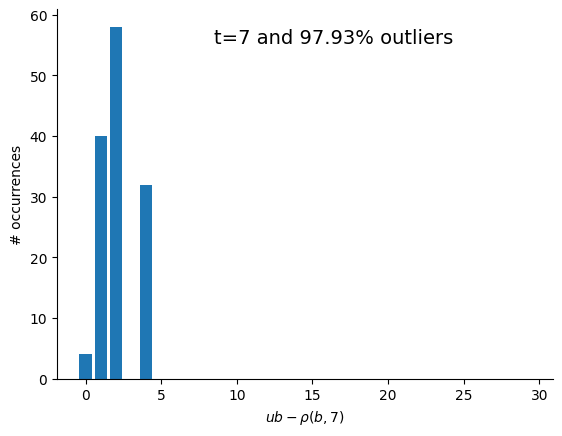} & \includegraphics[width=4cm]{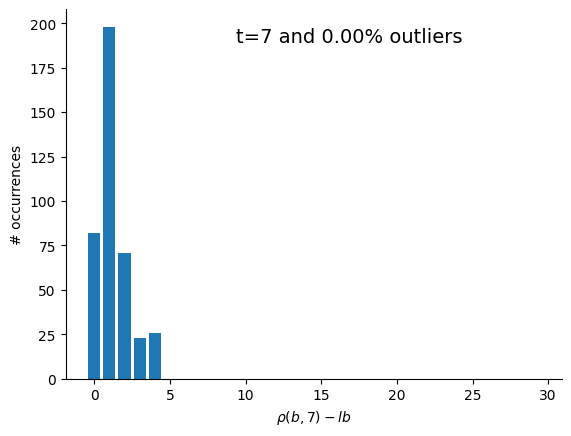} & \includegraphics[width=4cm]{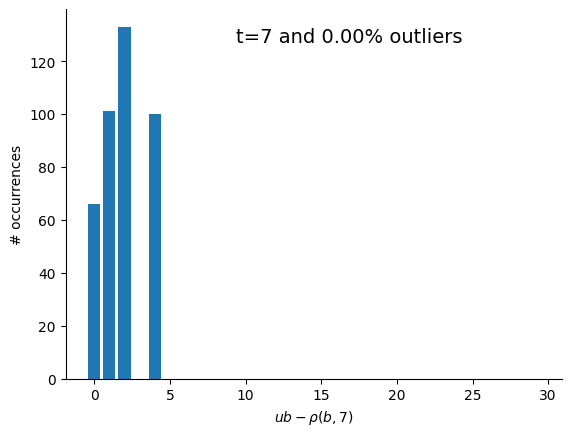} \\
  \end{tabular}
\end{center}
  \caption{Distribution of gap between extreme $\rho(b,t)$ and the lower (left)  and upper (right) bounds from Corollaries~\ref{elgamal_runs2} and~\ref{elgamal_runs1}, respectively, over all trials (two left columns) and with those where $g=v$ (two right columns). \label{runs_bound_dist}}
\end{figure}

We now turn our attention to runs. We give data on the behaviour of lower bounds from Corollary~\ref{elgamal_runs2} and upper bounds from Corollary~\ref{elgamal_runs1} since these are the better bounds for a majority of cases. Figure~\ref{run_bound_accuracy} shows the percentage of trials where 
there exist $b \in \z_v$ such that $\rho(b,t)$ matches the lower and upper 
bounds from Corollary~\ref{elgamal_runs2} and Corollary~\ref{elgamal_runs1} respectively. The blue bars show the numbers over 
all the trials. The green bars show the data for trials when $g=v$. Restricting 
to the trials where $\rho(b,t)>0$ is not as meaningful for runs and the accuracy 
for these cases match the lower bound accuracy over all trials very closely, 
with some loss of accuracy when $t=6,7$. Figure~\ref{runs_bound_dist} shows 
the distribution of $\min\{\rho(b,t):b \in \z_v\}-lb$ from Corollary~\ref{elgamal_runs2}and 
$ub - \max\{\rho(b,t):b \in \z_v\}$ from Corollary~\ref{elgamal_runs1} when these values are small over both all trials and just when $g=v$. The percentage of 
trials not appearing in the data is shown as the outliers in each case.  We 
ran experiments for the improvements given in Theorem~\ref{g_invertible} which 
resulted in essentially the same behaviour shown in Figure~\ref{runs_bound_dist}. 
The displayed data show that the bounds are frequently tight when $g=v$ and 
tight in at least some of trials when $g\neq v$.  As in the case for tuples, 
when the bounds are not tight, there are frequently not far off.  

\begin{figure}
  \includegraphics[width=8cm]{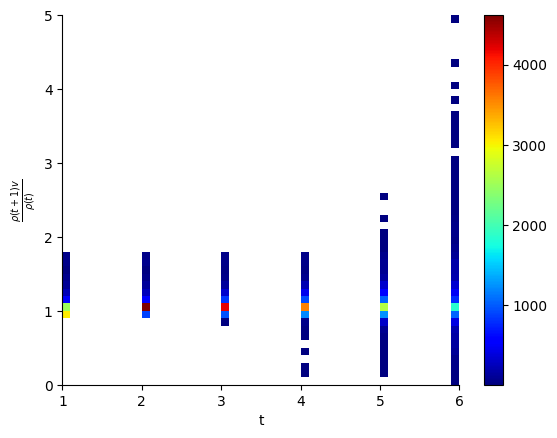} \includegraphics[width=8cm]{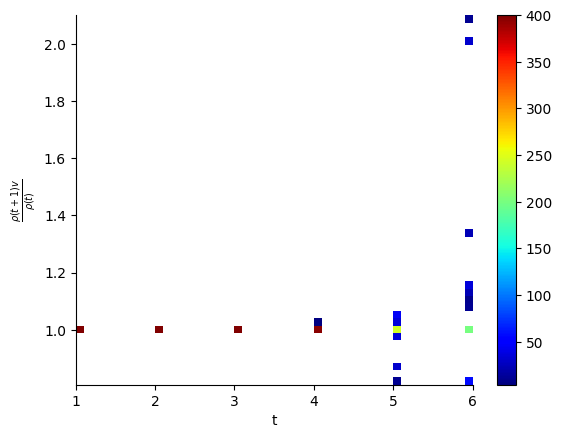}
  \caption{Distribution of values $\rho(t+1)v/\rho(t)$ over all trials (left) and when $g=v$ (right) shown as a heatmap.  Warmer colours indicate higher density. \label{ratios}}
\end{figure}

Finally we examine Golomb's run postulate experimentally. Figure~\ref{ratios} 
shows as heatmaps the ratios $\rho(t+1)v/\rho(t)$ over all the trials and over 
the trials where $g=v$. These show that ElGamal sequences modulo $v$ meet 
Golomb's run postulate quite well.

%%%%%%%%%%%%%%%%%%%%%%%%%%%%%%%%%%%%%%%%%%%%%%%%%%%%%%%%%%%%%%%%%%%%%%%%%%%%%
% \section{Applications}
% 
% \todo[inline]{Can we use this sequence as a pseudo-random bit generator? If so, 
% how good is it as a one time pad (counter mode block cipher). How hard is it
% to predict following term of the sequence?}

\section{Conclusions, discussion and further research} \label{conclusions}

%%%%%%%%%%%%%%%%%%%%%%%%%%%%%%%%%%%%%%%%%%%%%%%%%%%%%%%%%%%%%%%%%%%%%%%%%%%%%
% \section{Discussion}

Motivated by cryptographic concerns, von zur Gathen et al.~\cite{elgamalStats} 
investigated what properties of the ElGamal function are shared with random 
permutations. Experiments suggested that the distribution of cycle-sizes in 
these permutations match random permutations very well. Additionally they 
prove that the graph of this function, $\{(x,g^x): x \in \z_{p-1}\}$ is 
equidistributed which is also true of random permutations \cite{elgamalStats}. 
Similarly, Drakakis and collaborators in a series of papers look for measures 
that can distinguish the ElGamal function (and other Costas functions) from 
random permutations. They find that some common nonlinearity measures do not 
distinguish the two but that there are nonlinearity criteria for which the 
ElGamal function and other Costas functions are more nonlinear that expected 
by random \cite{drakakis_apn_2009,MR2723672,MR2660616}. We have continued 
these investigations by applying the remainder functions modulo $v$ to random 
permutations and the ElGamal function over $\z_p$, transforming them into 
sequences of length $p-1$ over $\z_v$. We then analyze both the random sequences 
and ElGamal sequences for their balance properties, period lengths, numbers 
of occurrences of tuples and runs which are all aspects of Golomb's 
randomness postulates \cite{GongGol2005}.

Proposition~\ref{balancedproperty} shows that every permutation of $\zpx$ 
produces a balanced sequence when reduced modulo $p$ and an exactly balanced sequences exactly when $p \equiv 1 \Mod v$. 
Thus, in this case, ElGamal sequences satisfy the balance randomness postulate. 
Purely random sequences can deviate from being balanced but the deviation 
is small for longer sequences. In Section~\ref{sec_balance_random} we have 
shown that the probability that a sequence $\sigma_v$ derived from a random 
permutation $\sigma$ has period less then $p-1$ is zero when 
$p \not\equiv 1 \Mod v$ and extremely small otherwise. In 
Theorem~\ref{thm_period_elgamal} we show that ElGamal sequences modulo $v$ 
always have period $p-1$. Thus, excepting for very unlikely events, the 
ElGamal sequences again behave the same as random balanced sequences.

In Theorem~\ref{thm_normal_random} we show that $\lambda(z)$ and $\rho(b,t)$ 
both have asymptotically normal distributions over the space of all random 
balanced sequence with period $p-1$. In Section~\ref{mean_variance} we compute 
the mean and standard deviation of these distributions. The leading term of 
both the expectation and variance of number of runs is $p/v^t$. Similarly the 
leading term of both the expectation and variance of $\rho(b,t)$ is $p/v^t$. 
Since $\lfloor x \rfloor = x + O(1)$, $\lceil x \rceil = x + O(1)$, 
Theorem~\ref{elgamal_tuple} shows that $\lambda(z)$ also has leading term 
$p/v^t$ for ElGamal sequences. Similarly Corollary~\ref{elgamal_runs2} shows 
that $\rho(b,t)$ has the same leading term for ElGamal sequences. Thus, for 
counts of both tuples and runs, the ElGamal sequences modulo $v$ match the 
average for random balanced periodic sequences.

Further, experiments for a range of values of $(p,g,v,t)$ show that
our lower and upper bounds on $\lambda(z)$ and $\rho(b,t)$ are tight
in at least some cases. When $g=v$ we prove that our lower bounds for
tuples are always tight and the cases where the upper bound is tight
can be precisely characterized.  When $g=v$ the bounds for runs are
significantly more accurate than when $g \neq v$. For both tuples and
runs, when the bounds are not tight the minimum and maximum values of
$\lambda(z)$ and $\rho(b,t)$ are frequently close to the bounds.  With
the current quantity of data the values of $\lambda(z)$ are not
normally distributed around their mean of $(p-1)/v^t$ but their
distribution does show some features consistent with such a
distribution.  This may be one aspect where ElGamal sequences differ
from sequences derived from random permutations. We think that are
larger set of experiments would clarify this question.  There is also
good evidence that Golomb's runs postulate is satisfied, especially
when $g=v$.

When $g$ is large with respect to $v$, the size of $D$ in the proof of  Theorem~\ref{elgamal_tuple} is quite large which means that there are many 
$X_d$. If the starting points, $\lceil pd_{t-1}/g^{t-1} \rceil$, of the 
various intervals of length $q+r/g^{t-1}$ from which  $x \equiv_v z(0)$ are 
drawn are evenly distributed modulo $v$, then the number of occurrences of 
$x$, either $\lfloor q/v \rfloor$ or $\lfloor q/v \rfloor +1$, even out 
when summed over $D$. In the case of runs of any symbol, we are summing 
over an even larger number of $X_d$ and different values of $z(0)$. 
Experiments encourage us that this ``smoothing out'' does occur for tuples 
and especially for runs but we have not been able to determine sufficient 
conditions. 

One promising scenario is if the sequence $\gamma_v$ begins with a run, that 
is $g^i\rem p \equiv_v g^{i+1}\rem p \neq g^0 = 1$, for $1 \leq i \leq t$. 
In this case the $x \in [1,p-1]$ that start runs are clustered so there will 
be $X_d$ which include ``runs'' of valid $x$ values which accomplishes some 
of the smoothing out we are after. Thus, we might be able to improve the 
floor or ceiling in the $q/v$ factor of the bounds. However the sizes of 
these clusters is not constant and this is a barrier to our efforts to prove 
improved bounds. Conversely it seems that if $\gamma$ does not begin with a 
run, which is the most frequent scenario,  then the starting points of runs 
are not clustered.  This means that we cannot use the clustering of $x \in X_d$. 
To improve the bounds in these normal cases we have to take advantage of 
the distributions of $\lceil pd_{t-1}/g^{t-1} \rceil$ over many $d$.

These are possible approaches to the broader problem of proving properties 
of the distributions of $\lambda(z)$ and $\rho(b,t)$ not just bounds on their 
range. The mean value for tuples is of course $(p-1)/v^t$, but what can we 
say theoretically about the variance of the distribution?  For runs, the mean 
and the variance should be studied. 

There are other important characteristics of a sequence including its
autocorrelation and linear complexity values that we do not include here. 
Our initial experiments on autocorrelation are mixed. Although the 
autocorrelation values are small in magnitude, there are a fairly large 
number of values obtained that may made them less attractive for applications
in communications. Our initial exploration of the linear complexity of our 
sequences shows that the linear complexity is high and very close to $(p-1)/2$. 
We leave for future work a complete study of these properties.

As we discuss throughout the article, $p \equiv 1 \Mod v$ is not necessary 
for many of our results. We state and prove all the results in 
Section~\ref{sec_elgamal} in full generality. In Section~\ref{sec_random} 
the results on period length and the fact of asymptotic normality of the 
distribution of counts of tuples and runs are fully general but we have 
only computed the means and variances when $p \equiv 1 \Mod v$. Similarly Lemma~\ref{many_to_one} only deals with the case that $v \mid (p-1)$. The 
latter can easily be generalized but the calculations from 
Section~\ref{mean_variance} would take a bit more effort.  We expect that 
the mean and variance will be slightly different but the leading term is 
likely the same and we leave this work for the future.

\section{Acknowledgments}

We would like to gratefully acknowledge the help, time and knowledge of 
Profs.~Jason Gao, Gennady Shaikhet and Yiqiang Zhao. Daniel Panario and 
Brett Stevens are supported by the Natural Sciences and Engineering Research 
Council of Canada (funding reference numbers RGPIN 05328 and 06392, 
respectively) and the Carleton-FAPESP SPRINT Program. Lucas Pandolfo Perin was supported by the Coordena\c{c}\~ao de 
Aperfei\c{c}oamento de Pessoal de N\'{i}vel Superior - Brasil (CAPES) - 
Finance Code 001.

\bibliographystyle{IEEEtran}
\bibliography{ref.bib}

\end{document}